\documentclass[reqno,a4paper,12pt]{article}

\usepackage{amsmath,amsthm,amsfonts,latexsym,amssymb,enumerate,color}
\usepackage{graphicx}
\usepackage{multirow}
\usepackage{float}
\usepackage[font=small,labelfont=bf]{caption}

\def\CC{\mathbb{C}}

\def\HH{\mathcal{H}}
\def\BB{\mathcal{B}}
\def\NN{\mathbb{N}}
\def\RR{\mathbb{R}}

\def\ZZ{\mathbb{Z}}

\def\Lp{L_{p}}
\def\Nev{{\rm Nev}}

\def\dim{\mathop{\rm dim}\nolimits}

\def\ker{\mathop{\rm ker}\nolimits}

\renewcommand\span{\mathop{\rm span}\nolimits}
\def\spam{\span} % I always mistype this

\def\im{{\rm Im}}
\def\ind{{\rm ind}}
\def\re{{\rm Re}}

\newtheorem{thm}{Theorem}[section]
\newtheorem{prop}[thm]{Proposition}
\newtheorem{defn}[thm]{Definition}
\newtheorem{lem}[thm]{Lemma}
\newtheorem{cor}[thm]{Corollary}

\newtheorem{ex}[thm]{Example}

\numberwithin{equation}{section}

\def\beginpf{\begin{proof}}
\def\endpf{\end{proof}}
\def\beq{\begin{equation}}
\def\eeq{\end{equation}}

\def\imag{\mathop{\rm Im}\nolimits}

\begin{document}

\title{Kernels of unbounded Toeplitz operators and factorization of symbols}

\date{\today}

\author{M.~C. C\^amara\thanks{
Center for Mathematical Analysis, Geometry, and Dynamical Systems,
Departamento de Matem\'{a}tica, Instituto Superior T\'ecnico, 1049-001 Lisboa, Portugal. \tt ccamara@math.ist.utl.pt},\quad
 M.T. Malheiro\thanks{Centre of Mathematics, Departamento de Matem\'atica, Universidade do Minho, Campus de Azur\'em, 4800-058 Guimar\~aes, Portugal. \tt mtm@math.uminho.pt},\quad and
 J. R. Partington\thanks{School of Mathematics, 
 University of Leeds, Leeds LS2~9JT, U.K. \tt j.r.partington@leeds.ac.uk}\ \thanks{Corresponding author}}

\maketitle

\begin{abstract}
We consider kernels of unbounded Toeplitz operators in $H^p(\mathbb C^+)$  in terms of a factorization of their symbols. We study the existence of a minimal Toeplitz kernel containing a given function in $H^p(\mathbb C^+)$, we describe the kernels of Toeplitz operators whose symbol possesses a certain factorization involving two different Hardy spaces and we establish relations between the kernels of two operators whose symbols differ by a factor which corresponds, in the unit circle, to a non integer power of $z$. We apply the results to describe the kernels of all Toeplitz operators with non-vanishing piecewise continuous symbols.

\end{abstract}



%%%%%%%%%%%%%%%%%%%%%%%%%%%%%%%%%%%%%%%%%%%%%%%%%%%%%%%%%%%% Introduction

\section{Introduction}
\label{sec:intro}

In   \cite{Sarason}, Sarason presented the basic theory of unbounded Toeplitz operators in $H^2(\mathbb{D})$ with symbols in $L^2(\mathbb{T})$ and, motivated by natural questions   that lead to other types of symbols \cite{Helson, Suarez, Seubert}, of Toeplitz operators with analytic and co-analytic symbols in more general classes.
Unbounded Toeplitz operators appear naturally, for instance when studying inverses or generalized inverses of Toeplitz operators with bounded symbols \cite{5GK, 8MP}. Indeed, the inverse of a bounded Toeplitz operator, if it exists, is the composition of two Toeplitz operators which, in general, are unbounded.

\vspace{3mm}

We can consider, analogously, Toeplitz operators in the Hardy space $H^2(\mathbb{C}^+)$ of the upper half-plane $\CC^+$, and more generally, in $H^p(\mathbb{C}^+)$, $p>1$, which arise in many applications  \cite{2Dud,6Kal,8MP,11S,12T,13V,14W}. One may ask, in that case, what are the natural classes of symbols to consider and what properties do those operators possess. In particular, we would like to examine their kernels and study what properties are shared with kernels of bounded Toeplitz operators, which have attracted great interest for their rich structure and the information that they provide on the corresponding Toeplitz operators (see for instance the recent survey paper \cite{HM}).

\vspace{3mm}

We assume here that $g$ is a measurable function %satisfying the condition
%\beq \label{Introd-1}
%\int_{\mathbb{R}}\frac{\log |g(t)|}{1+t^2}dt < \infty\eeq
 and that there exists a non-zero $f \in H^p(\mathbb{C}^+)=H_p^+$ such that $gf \in L^p(\mathbb{R})=L_p$. We define the Toeplitz operator $T_g$ on the domain
 \[D(T_g)=\left\{\phi_+ \in H_p^+: g \phi_+ \in L_p\right\}\]
by
\[T_g \phi_+= P^+g\phi_+, \qquad \text{ for } \phi_+ \in D(T_g),\]
where $P^+$ is the Riesz projection from $L^p$ onto $H_p^+$.
Then $\ker T_g$ consists of the functions $\phi_+ \in H_p^+$ such that
\[ g\phi_+=\phi_- \qquad \text{ with } \phi_- \in H_p^-=\overline{H_p^+}.\]

The kernels of Toeplitz operators are also called Toeplitz kernels.
Given any function $\phi_+ \in H_p^+$, there is always a Toeplitz kernel containing $\phi_+$. In the class of bounded Toeplitz operators, for each
$\phi_+ \in H_p^+ \setminus \{0\}$ one can find a minimal kernel $\mathcal{K}_{\min}(\phi_+)$, which is contained in any other Toeplitz kernel to which $\phi_+$
belongs. If $\phi_+=I_+O_+$ is an inner-outer factorization of $\phi_+$, with $I_+$ inner and $O_+$ outer in $H_p^+$, then $\mathcal{K}_{\min}(\phi_+)= \ker T_{\overline{I_+}{\overline{O_+}}/{O_+}}$ \cite{CP14}. A first question that can be asked when we consider unbounded Toeplitz operators is whether we can also find a minimal kernel containing $\phi_+$, for general symbols, or whether one can find classes of possibly unbounded symbols
% satisfying  \eqref{Introd-1},
for which the minimal kernel of $\phi_+$ exists and can be compared with $\mathcal{K}_{\min}(\phi_+)$, for instance by inclusion.

\vspace{3mm}

A second question is how to determine whether or not $\ker T_g$ is trivial, which is equivalent to the question of injectivity of $T_g$. The nontriviality of the kernel is directly connected with the existence of certain types of factorization of the symbol. Indeed the existence of a non-zero function $\phi_+$ in  $\ker T_g$
means that $g\phi_+=\phi_- \in H_p^-\setminus \{0\}$. If $\phi_+=I_+O_+$ as above, and $\phi_-=I_-O_-$ with $\overline{I_-}$ inner and  $\overline{O_-}$ outer in $H_p^+$, we must then have
\beq \label{Introd-5} g=O_-\overline{I}O_+^{-1}\eeq
for some inner function $I$ and $O_+$, $\overline{O_-}$ outer in $H_p^+$, and conversely, if $g$ admits a factorization \eqref{Introd-5} then $I_+O_+ \in \ker T_g$ for
each $I_+$ that divides $I$. It is clear that if $I$ in \eqref{Introd-5} is not a finite Blaschke product, then $\dim \ker T_g=\infty$.
We will  thus consider symbols possessing some factorization of the form \eqref{Introd-5}. It is very difficult, however, to describe the kernels of operators whose symbol admits such a general representation, without imposing certain conditions on the inverses of the factors.

For bounded symbols, it is well known that we can go further and study their invertibility and Fredholmness based on an appropriate factorization of their symbols, where conditions are imposed on the factors as well as their inverses (\cite {C,2Dud,5GK}).   In order to define this factorization, let $\mathcal{R}$ denote the set of rational functions belonging to $L^{\infty}(\mathbb{R})=L_{\infty}$,
and write
\[
\mathcal{R}_0=\mathcal{R}\cap L_p, \quad \lambda_{\pm}(\xi)=\xi \pm i, \quad r(\xi)=\frac{\xi-i}{\xi+i}  \hbox{ for }\xi \in \mathbb{R}, \quad \hbox{and }
 \mathcal{H}_p^{\pm}=\lambda_{\pm}H_p^{\pm}.
 \]

Note that  \eqref{Introd-5} can be rewritten as a product of the form $g=g_-\overline{Ir}g_+$ with $g_- \in \mathcal{H}_p^-$, $g_+^{-1} \in \mathcal{H}_p^+$.

%Definition 1.1
\begin{defn} \label{defn_int_1} By a $p$-factorization of $g\in L_{\infty}$ we mean a representation of $g$ as a product
\beq \label{Introd-7} g=g_-r^kg_+\eeq
where $k \in \mathbb{Z}$ and
\[g_+ \in \mathcal{H}_{p'}^+, \quad g_+^{-1} \in \mathcal{H}_{p}^+, \quad g_- \in \mathcal{H}_{p}^-, \quad  g_-^{-1} \in \mathcal{H}_{p'}^-\]
with $1/p+1/p'=1$; if moreover
\beq \label{Introd-9} g_+^{-1}P^+g_-^{-1}I: \mathcal{R}_0\rightarrow L_p \quad \text{ is bounded, }\eeq
where $I$ denotes the identity operator,
then \eqref{Introd-7} is called a Wiener--Hopf $p$-factoriz\-ation. The integer $k$ is called the index of the factorization, and if $k=0$
the factorization is said to be canonical.
\end{defn}

 If $g$ admits a $p$-factorization then it is unique, apart from constant factors. The existence of a  $p$-factorization (where condition
\eqref{Introd-9} is not taken into account) allows one
to characterize the kernel of $T_g$ and that of $T_g^*$; the operator $T_g$, $g \in L_{\infty}$, is Fredholm if and only if $g$ admits a Wiener--Hopf $p$-factorization and, in that case, its Fredholm index is $-k$; the operator is invertible if and only if $g$ admits a canonical Wiener--Hopf $p$-factorization  \cite{C,LS}.

%\vspace{3mm}

A $p$-factorization  may not exist, however, when $g$ is not invertible in the algebra of functions continuous on $\mathbb{R}_{\infty}=\mathbb{R}\cup\{\infty\}$. This is the case of symbols as simple as $r^{1/2}(\xi)=\sqrt{\frac{\xi-i}{\xi+i}}$, where we assume a discontinuity at $\infty$, for $p=2$; it is easy to see that $r^{1/2}$ does not admit a 2-factorization. We thus generalize Definition \ref{defn_int_1} and study symbols possessing a $(j,s)$-factorization, defined in Section \ref{sec_pq_fact}, which coincides with a $p$-factorization if $j=p\,,\,s=p'$ with $1/p+1/p'=1$. This will allow us to describe the kernels of a wide class of Toeplitz operators, including unbounded symbols and piecewise continuous symbols that do not admit a $p$-factorization, and establish criteria for a Toeplitz kernel to be trivial or not.

\vspace{3mm}

A third question regarding the kernels of Toeplitz operators with a possibly unbounded symbol is the relation between $\ker T_g$ and $\ker T_{r^{\alpha}g}$, where $\alpha \in \mathbb{R}$.  The relation between $\ker T_g$ and $\ker T_{r^{k}g}$, with $k \in \mathbb{Z}$, was studied in  \cite{BCD, CMP} for bounded symbols. It makes sense to ask also what happens if we multiply a symbol by some non-integer power of $r$. Indeed taking for instance the algebra of all piecewise continuous functions in $\mathbb{R}_{\infty}$, every function in that class can be represented as the product of a continuous function, whose kernel can be described from a $p$-factorization, by as many non-integer powers of $r$ as the number of points of discontinuity in $\mathbb{R}$ (\cite{2Dud}).

\vspace{3mm}

In this paper we address the three questions mentioned above. In Section \ref{sec:2} we begin by recalling some
preliminary results on Wiener--Hopf factorizations. Then in Section \ref{sec:3} we present some aspects of the theory of unbounded
Toeplitz operators, with particular reference to domains and kernels.
Section \ref{sec:4} is concerned with developing the theory of minimal kernels.
In Section \ref{sec:5} we introduce the notion of a factorization
based on two indices, giving a generalization of $p$-factorization, and  Section \ref{sec:6}
studies the properties of the corresponding factors.
Next, in Section \ref{sec:7} we apply the theory of factorization to study the kernels of  Toeplitz operators with both unbounded and bounded symbols. Finally, Section \ref{sec:8} addresses
questions regarding the relation between $\ker T_g$ and $\ker T_{r^{\alpha}g}$, where $\alpha \in \mathbb{R}$, obtaining a description of the kernels of Toeplitz operators in $H^p$ with  piecewise continuous symbols that do not possess a $p$-factorization.

% Section 2 - Preliminaries

\section{Preliminaries}
\label{sec:2}

Let $\mathcal{GA}$ denote the group of invertible elements in an algebra $\mathcal{A}$ and let $C(\mathbb{R}_\infty)$ represent the algebra of all functions which are continuous in $\mathbb{R}_{\infty}=\mathbb{R} \cup \{ \infty \}$ and possess equal limits at $\pm \infty$.

For $1\leq p<\infty$, let $\Lp$ denote the Lebesgue space of all complex Lebesgue measurable functions $f$ which are $p$-integrable in $\RR$ %, with the norm
%\[
%\|f\|_{p}=\left(\int_{\RR}|f(x)|^p\, dx\right)^{1/p}.
%\]
%Consider $p'$ such that $\frac{1}{p}+\frac{1}{p'}=1$.
%We denote
and by $H_p^+$ the Hardy space $H^p(\CC^+)$ of the upper half plane $\imag (z)<0\,,\,z \in \mathbb C$. %all functions $f$ which are analytic in the upper half plane $\CC^+$ and such that, for all $y>0$, $|f(x+iy)|^p$ is integrable over $\RR $ and there is a constant $M \in \RR^+$ such that
%\[
%\int_{\RR}|f(x+iy)|^p dx <M \text{ for all } y>0.
%\]
%$H_p^+$ is a Banach space with the norm
%\[
%\|f\|_{H_p^+}=\sup_{y\in \RR^+}\left(\int_{\RR}|f(x+iy)|^p\, dx\right)^{1/p}.
%\]
 We define $H_p^-$ similarly for the lower half plane $\CC^-$ and we identify as usual  $H_p^\pm$ with closed subspaces of $L_p$.  %Identifying each function $f_{\pm} \in H_p^{\pm}$ with its nontangential boundary value on $\RR$, %\cite{3Dur},
%we have that $H_p^+$ and $H_p^-$ are closed subspaces of $L_p$ and $H_p^-=\overline{H_p^+}$. For $p>1$,
%\[H_p^+\cap H_p^-=\{0\}\,\,,\,\, L_p=H_p^- \oplus H_p^+,\]
%where the sum is orthogonal if $p=2$.
We denote by $P^\pm$ the projections from $L^p$ onto $H_p^\pm$, parallel to $H_p^\mp$. %We have
%\[P^++P^-=I, \qquad P^+-P^-=S_{\mathbb{R}},\]
%where $\dps S_{\mathbb{R}}f(\xi)=\frac{1}{\pi i}\int_{\mathbb{R}}\frac{f(t)}{t-\xi}dt$ (the latter being understood as a Cauchy principal value).

It is well known that for $p, r \in ]1,+\infty[$, if $\,
f \in H_p^{\pm}$ and $g \in H_r^{\pm}$ then $fg \in  H_s^{\pm}$ with $ \frac{1}{p}+\frac{1}{r}=\frac{1}{s}.$
Recall that
$\lambda_\pm(\xi)=\xi\pm i\,,\, r(\xi)=\frac{\xi+i}{\xi-i}.$
 We define, for $p \geq 1$,
\[\mathcal{L}_p=\lambda_+L_p, \qquad \mathcal{H}_p^{\pm}=\lambda_\pm H_p^{\pm},\]
\[\mathcal{B}_p=\lambda_+^2L_p, \qquad \mathcal{B}_p^{\pm}=(\lambda_\pm)^2H_p^{\pm}.\]
Then we have (\cite{14W})
\begin {equation}\label {2.3a}
\HH_p^+ \cap \HH_p^-=\CC \quad \text{if } p>1; \qquad \HH_p^+ \cap \HH_p^-=\{0\} \quad \text{if } p=1.
\end {equation}
If $f\in \BB_p^+ \cap \BB_p^-$ then there exist $f_\pm \in \HH_p^\pm$ such that
$\frac{\xi+i}{\xi-i}f_+=f_-$, therefore we have
$\frac{\xi+i}{\xi-i}\,(f_+- f_+(i))=f_- -\frac{f_+(i)}{\xi-i}$,
where the left-hand side belongs to $\HH_p^+$  and the right-hand side belongs to $\HH_p^-$. It follows from \eqref {2.3a} that both sides are equal to a constant which is zero if $p=1$. Consequently,

\beq \label{BB_inters}\BB_p^+ \cap \BB_p^-=\mathcal{P}_1 \quad \text{if } p>1; \qquad \BB_p^+ \cap \BB_p^-=\CC \quad\text{if } p=1,\eeq
where $\mathcal{P}_1$ represents the space of polynomials with degree less or equal to 1.

We denote by $L_{\infty}$ the space of all essentially bounded functions on $\RR$
%, with the norm
%\[\|f\|_{\infty}= \ess \sup_{x \in \RR} |f(x)|\]
and by $H_{\infty}^+ \, (H_{\infty}^-)$ the space of all functions analytic and bounded in $\mathbb{C}^+ (\mathbb{C}^-)$. We identify $H_{\infty}^{\pm} $ with the subspaces of $L_{\infty}(\mathbb{R})$ consisting of their (nontangential)  boundary functions on $\RR$.

%Let us call $r$ the function defined on $\RR$ by $r(\xi)=\frac{\xi-i}{\xi+i}$ and $r_{\pm}$ the functions defined by
%\begin{equation}
%r_{\pm}(\xi)=\frac{1}{\xi\pm i}, \quad \xi\in \mathbb{R}.
%\end{equation}
%Let $P$ and $Q$ be the complementary projections on $\Lp$ such that, for $f \in H_p^+$, \beq \label{2.8}
%Pf=f, \qquad Qf=0.
%\eeq

%Analogously, if $f \in H_p^-$, then
%\begin{equation}\label{2.9}
%Qf=f, \qquad Pf=0.
%\end{equation}
%It is also well-known that $P$ and $Q$ defined as in \eqref{2.8} and \eqref{2.9}, are continuous projections mapping $\Lpr$, $p>1$, onto $H_p^+$, $H_p^-$, respectively, such that

For $p>1$, let $p'$ be defined by $\frac{1}{p}+\frac{1}{p'}=1$.
In what follows we will always assume that $p \in ]1,+\infty[.$

%Def 2.1
\begin{defn} By a \emph{$p$-factorization} of $g \in L_{\infty}$ we mean a representation of $g$ as a product
\beq \label{def_p_fact}
g=g_-r^kg_+
\eeq
where $k \in \mathbb{Z}$ and
\beq
g_+ \in \HH^+_{p'}, \, g_+^{-1} \in \HH^+_{p}, \, g_- \in \HH^-_{p}, \, g_-^{-1} \in \HH^-_{p'}.
\eeq
\end{defn}
If moreover,
\beq \label{inv_op}
g_+^{-1}P^+g_-^{-1}I: \mathcal{R}_0\rightarrow L_p  \quad \text{is bounded,}
\eeq
 where $\mathcal{R}_0$ denotes the space of all rational functions in $C(\mathbb{R}_{\infty})\cap L_p$ then \eqref{def_p_fact} is called a \emph{Wiener--Hopf (WH) $p$-factorization}, or a \emph{generalized factorization relative to $L_p$}. The integer $k$ in \eqref{def_p_fact} is called the factorization index. If $k=0$ then the factorization is said to be \emph{canonical}.

The representation \eqref{def_p_fact}, if it exists, is unique up to non-zero constant factors. It is well known \cite{C, CDR, 2Dud, 5GK} that the Toeplitz operator with symbol $g \in L_{\infty}$,
\beq
T_g : H_p^+ \rightarrow H_p^+, \quad T_g(\phi_+) = P^+g \phi_+
\eeq
is Fredholm if and only if $g$ admits a WH $p$-factorization; it is invertible if and only if the WH $p$-factorization is canonical and, in that case, $T_g^{-1}$ is the bounded extension of the operator \eqref{inv_op} to $H_p^+$:
\[T_g^{-1}=g_+^{-1}P^+ g_-^{-1}: H_p^+ \rightarrow H_p^+.\]
If $g$ admits a $p$-factorization, not necessarily satisfying \eqref{inv_op}, then this is enough to describe the kernel of the Toeplitz operators $T_g$ and its adjoint $T_g^*=T_{\overline{g}}$. We have
\beq
\ker T_g = 0 \text{ if } k \geq 0, \quad \ker T_g= g_+^{-1}K_{r^{|k|}} \text{ if } k<0
\eeq
where by $K_{r^{n}}$, $n \in \mathbb{N}$, we denote the model space
\beq
K_{r^{n}}=\spam\left\{\frac{1}{\xi+i}\left(\frac{\xi-i}{\xi+i}\right)^j, j=0, \ldots, n-1 \right\}.
\eeq
It may happen that two factorizations of the same function with respect to different $L_p$ spaces coincide, as it happens with continuous functions. We have the following

%Theor 2.2
\begin{thm} If $g \in \mathcal{G}C(\RR_{\infty})$ then $g$ admits a WH $p$-factorization, which is the same for any $p>1$. The factorization index is equal to the winding number of $g$ with respect to the origin.
\end{thm}

A $p$-factorization may not exist, however, for bounded functions as simple as a non-integer power of $r$.

For $c \in \mathbb{R}_{\infty}$, $\alpha \in \mathbb{C} \setminus \mathbb{Z}$, we define
\[z_c^{\alpha}=e^{\alpha \log z}\]
where the branch cut connecting 0 to $\infty$ intersects the unit circle at the point $\displaystyle z_0=\frac{c-i}{c+i}$, with $z_0=1$ if $c=\infty$, and we take $\log 1=0$. Then the function
\beq \label{r_form}
r^{\alpha}_c(\xi)=\left(\frac{\xi-i}{\xi+i}\right)_c^{\alpha}
\eeq
is continuous for all points of $\mathbb{R}_{\infty}$ except for the point $c$ where it has different finite one-sided limits. Other piecewise continuous functions can be expressed as products of a continuous function by non-integer powers of $r$ of the form \eqref{r_form}. Let $PC(\mathbb{R}_{\infty})$ denote the space of all piecewise continuous functions $g$, with finite limits at $\pm \infty$. Any $ g \in PC(\mathbb{R}_{\infty})$ that does not vanish on $\mathbb{R}_{\infty}$ can be represented as a product
\beq \label{g_decomp}
g = h  \, r_{\infty}^{\alpha_{\infty}}\prod_{j=1}^n r_{c_j}^{\alpha_j}\eeq
where $h \in \mathcal{G}C(\mathbb{R}_{\infty})$, $c_j \in \mathbb{R}$ ($j=1, \ldots, n$) are the points of discontinuity of $g$ on $\mathbb{R}$ and
\beq \label{alpha_inf_decomp}
\alpha_{\infty}=\frac{1}{2 \pi i}\log \frac{g(+\infty)}{g(-\infty)} \quad \text{ with } \, -\frac{1}{p}< \re \, \alpha_{\infty} \leq \frac{1}{p'} \, ,
\eeq
\beq\label{alpha_j_decomp}
\alpha_{j}=\frac{1}{2 \pi i}\log \frac{g(c_j^-)}{g(c_j^+)} \quad \text{ with } \, -\frac{1}{p'}< \re \, \alpha_{j} \leq \frac{1}{p}, \, (j=1, \ldots, n).
\eeq
(\cite{2Dud}, see also \cite{C}). Without loss of generality, one may assume that $\alpha_{\infty}$ and $\alpha_{j}$, $j=1, \ldots, n,$ are real. Indeed $r^{i \im (\alpha_{j})}$ is real and positive and therefore admits a bounded factorization that can be absorbed into $h$. Denoting by $m$ the number of exponents $\alpha_{j}$ in \eqref{g_decomp} that satisfy the condition $-\frac{1}{p'}< \re\, \alpha_{j} < \frac{1}{p}$, we can also write \eqref{g_decomp} as
\beq \label{g_decomp_b}
g = h \, r_{\infty}^{\alpha_{\infty}}\prod_{j=1}^m r_{c_j}^{\alpha_j}\prod_{j=1}^s r_{c_j}^{1/p}\,.\eeq
To each $g \in PC(\mathbb{R}_{\infty})$ and each $p>1$ we associate the function \cite{2Dud}
\beq \label{g_p}
g_p(\xi, w)=\frac{g(\xi^-)+g(\xi^+)}{2}+\frac{g(\xi^-)-g(\xi^+)}{2}\coth(\pi(\frac{i}{p}+w)),
\eeq
$\xi \in \RR_{\infty}, w \in \RR \cup \{\pm \infty\}=\RR_{\pm \infty}$. Here we take $g(\infty^{\pm})=g(\mp \infty)$. If $g \in C(\RR_{\infty})$ then $g_p(\xi, w)=g(\xi)$ for all $\xi \in \RR_{\infty}$, $w \in \RR_{\pm \infty}$, $p>1$. It is easy to see that if $p=2$ then \eqref{g_p} takes the form
\beq \label{g_2}
g_2(\xi, w)=\frac{g(\xi^-)+g(\xi^+)}{2}+\frac{g(\xi^-)-g(\xi^+)}{2}\tanh(\pi w).
\eeq
The image of $g_p$ in the complex plane is a closed curve $\Gamma$ obtained by adding to the image of $g(\xi)$, with $\xi \in \RR_{\infty}$, certain arcs of a circle (or line segments if $p=2$) connecting the points $g(\xi^-)$ and $g(\xi^+)$ whenever these two values are different. If
\[\inf_{(\xi, w) \in \RR_{\infty} \times \RR_{\pm \infty}} |g_p(\xi,w)|>0,\]
which means that the image of $g_p$ in the complex plane, $\Gamma$, is a closed curve that does not pass by zero, then we say that $g$ is \emph{$p$-nonsingular}. In this case, we associate with $g_p$ an integer, $\ind \, g_p$, which is the winding number of $\Gamma$ around the origin.

If two functions $f$ and $g$ belonging to $PC(\RR_{\infty})$ have no common points of discontinuity and are $p$-nonsingular, then $fg$ is also non-singular and $\ind_p (fg)=\ind_p f+\ind_p g$ (\cite{2Dud}).

We have the following:
%Theor 2.3
\begin{thm} \cite{2Dud} Let $g \in PC(\RR_{\infty})$ of the form  \eqref{g_decomp} - \eqref{alpha_j_decomp}. Then $g$ is $p$-nonsingular if and only if $-\frac{1}{p}<\re\, \alpha_{\infty}<\frac{1}{p'}$, $-\frac{1}{p'}<\re\, \alpha_{j}<\frac{1}{p}$ for all $j=1, \ldots, n$.
\end{thm}

%THeor 2.4
\begin{thm} \cite{2Dud} Let $g \in PC(\RR_{\infty})$. The operator $T_g$ has closed range in $H_p^+$ if and only $g$ is $p$-nonsingular. In that case $T_g$ is Fredholm, with Fredholm index $-\ind_p g$, invertible if $\ind_p g=0$, and
\[g=\left[h_- (\xi-i)^{\alpha_{\infty}}\prod_{j=1}^n\left(\frac{\xi-i}{\xi-c_j}\right)^{\alpha_j}\right]r^k\left[h_+ \frac{1}{(\xi+i)^{\alpha_{\infty}}}\prod_{j=1}^n\left(\frac{\xi-c_j}{\xi+i}\right)^{\alpha_j}\right]\]
is a $p$-factorization of $g$.
\end{thm}

%Section 3   Unbounded Toeplitz operators

\section{Unbounded Toeplitz operators}
\label{sec:3}

Let $g$  be a measurable function on $\RR$
such that there exists $f_+ \in H_p^+ \setminus \{0\}$ with $gf_+ \in L_p$. We denote this class of functions by $\sigma_p$. Let $T_g$ be the (possibly unbounded) Toeplitz operator  defined on the domain
\beq \label{1}
D(T_g)=\{f_+ \in H_p^+: gf_+ \in L_p\}
\eeq
by \beq \label{1tg}T_gf_+=P^+(gf_+).\eeq

Note that, if we assume $g$ to be a measurable function on $\mathbb{R}$ without any further restriction, it is possible for $D(T_g)$ to be $\{0\}$. For example,
if $g(\xi)=e^{\xi^2}$ there is no function $f_+ \in H_p^+ \setminus \{0\}$ with $gf_+ \in L_p$. In fact,

\beq \label{2}
\int_{\RR} \frac{|\log W(t)|}{1+t^2}\,dt < \infty
\eeq
 is a necessary and sufficient
condition for there to be an outer function $O$ such that $W=|O|$ (\cite{Nikolski}, sections 3.9 and 6.4); explicitly, to within a constant, the outer function is
\beq \label{3}
O(\xi)=\exp\left(\frac{1}{i\,\pi}\int_{\RR}\frac{\log W(t)}{1+t^2}\frac{1+t\xi}{t-\xi}\,dt\right).
\eeq

So, for  $f_+ \in H_p^+\setminus \{0\}$, we can have $\frac{\log |f_+(\xi)|}{1+ \xi^2} \leq -1$
only on a set of finite measure, so that $|f_+(\xi)|\geq \exp(-1-\xi^2)$, i.e.,  $|e^{\xi^2}f_+(\xi)|\geq e^{-1}$ except on a set of finite measure and therefore $e^{\xi^2}f_+$ cannot lie in $L_p$.

%Proposition 3.1
\begin{prop} $T_g=0$ if and only if $g=0$.
\end{prop}

\beginpf
Assume that  $T_g=0$. Then there exists $f_+ \in H_p^+ \setminus \{0\}$ such that $gf_+ \in L_p$, so $r^nf_+ \in D(T_g)$ for all $n \in \NN$, and we have $gf_+ \in H_p^-,$ $r^ngf_+ \in H_p^-$ for all $n \in \NN$, implying that $gf_+=0$. Since
$f_+ \in H_p^+ \setminus \{0\}$, it follows that $g(\xi)=0$ a.e..
\endpf

From now on we assume that $g \neq 0$ and $g\in \sigma_p$.

The kernel of $T_g$ is the subspace
\beq \label{tgkernel}
\ker T_g =\left\{f_+ \in H_p^+ : gf_+ \in H_p^-\right\}.
\eeq
It is clear that $\ker T_g$ is nearly $S^*$-invariant, i.e., a subspace $M$ of $H_p^+$ such that
\beq
f \in M, \quad r^{-1}f \in H_p^+ \Rightarrow r^{-1}f \in M.
\eeq

Condition \eqref{2} provides a necessary condition for this kernel to be non-zero. Indeed, if $f_+ \in H_p^+ \setminus \{0\}$ and $gf_+ \in H_p^-$, then $gf_+\in H_p^-\setminus\{0\}$ and \eqref{2} is satisfied for both $W=|f_+|$ and $W=|gf_+|$; so we automatically have the same for $W=|g|$. We have then the following.
%Proposition 3.2
\begin{prop}
A necessary condition for $\ker T_g$ to be different from $\{0\}$ is that
\beq \label{1cond}
\int_{\RR} \frac{|\log |g(t)||}{1+t^2}\,dt < \infty.
\eeq
\end{prop}

%Clearly, if $g \in (\xi+i)^m L_r$ for some $m \in \NN$, $r>0$, then $g$ belongs to the class of symbols satisfying \eqref{1cond}.

%Proposition 3.3
\begin{prop}
If $g$ satisfies \eqref{1cond} and $O$ is the outer function defined by \eqref{3}, with $W=|g|=|O|$, then
\[D(T_g)=\{f_+ \in H_p^+ : Of_+ \in H_p^+\}\]
\end{prop}

\beginpf Clearly, $\{f_+ \in H_p^+ : Of_+ \in H_p^+\}\subset D(T_g)$. To show the converse, note that $O$ is the ratio of two $H_p^+$ outer functions, obtained by taking $W=|gf_+|$ and $W=|f_+|$ in \eqref{3}. Thus $O$ belongs
to the Smirnov class $\Nev_+$ \cite{Nikolski} and if $gf_+ \in L_p$ then $Of_+ \in L_p \cap \Nev_+=H_p^+$.
\endpf

%Proposition 3.4
\begin{prop}\label{prop3} \begin{sloppypar} If $g$ satisfies \eqref{1cond} then there exists an outer function $Q$, bounded below,  such that  $f_+ \mapsto Qf_+$  is an isometry from $D(T_g)$, endowed with the norm $\|f_+\|_{Q}=\|Qf_+\|_{p}$, onto $H_p^+$. Moreover, $T_g$ is bounded on $\left(D(T_g), \|\cdot \|_{Q}\right).$ \end{sloppypar}
\end{prop}
 \beginpf The result holds provided that $f_+\in D(T_g)\Leftrightarrow Qf_+ \in H_p^+$. Suppose that there is an outer  function $O$, bounded below, with $|O|=|g|$. Then $f_+ \in D(T_g) \Leftrightarrow Of_+ \in H_p^+$ and the result holds with $Q=O$. If $O$ is not bounded below we can choose an outer function $Q$ such that $|Q|=1+|O|=1+|g|$, which is bounded below. If $f_+ \in D(T_g)$ then $f_+ \in H_p^+$ and  $Of_+ \in H_p^+$,  so $Qf_+ \in H_p^+$ because $Qf_+ \in L_p$ and $Qf_+ =f_+/Q^{-1}$ (where $Q^{-1}$ is outer in $H_{\infty}^+$) is in $\Nev_+$. Conversely, since $Q$ is bounded below, if $Qf_+ \in H_p^+$ then $f_+ \in H_p^+$. Therefore we have that $f_+ \in D(T_g) \Leftrightarrow Qf_+ \in H_p^+$ and the result holds.
 \endpf

For $g \in L_{\infty}$, $\ker T_g$ is a nearly $S^*$-invariant closed subspace of $H_p^+$ \cite{CP14}. For more general symbols, we have the following.

%Corollary 3.5
\begin{cor} With the same assumptions as in Proposition \ref{prop3}, $Q \ker T_g$ is a closed nearly $S^*$-invariant subspace of $H_p^+$ and, for $p=2$, $\ker T_g=Q^{-1}FK_{\theta}$  where $\theta$ is an inner function, $K_{\theta}=\ker T_{\bar{\theta}}$ is the model space defined by $\theta$, $F$ is an isometric (outer) multiplier from $K_{\theta}$ onto $Q \ker T_g$.
\end{cor}

\beginpf $M=\ker T_g$ is a $\|\cdot\|_Q$  closed and nearly $S^*$ - invariant subspace of $D(T_g)$, so $QM$ is a closed nearly $S^*$ - invariant subspace of $H_p^+$. In the case $p=2$ the result follows from the Hayashi--Hitt results \cite{hayashi,hitt}.
\endpf

%Section 4
\section{Minimal kernels}
\label{sec:4}

Given any $f_+ \in H_p^+\setminus \{0\}$, there is always a Toeplitz kernel containing $f_+$; moreover there exists a minimal kernel
that contains $f_+$ and is contained in any other Toeplitz kernel, with bounded symbol, to which $f_+$ belongs. It is denoted by $K_{\min}(f_+)$ and we can associate with it a unimodular symbol: %$g=\overline{I_+}\overline{O_+}/O_+$:
\beq \label{4.1aa}
K_{\min}(f_+)=\ker T_{\overline{I_+}\overline{O_+}/O_+}
\eeq
(\cite{CP14}). The function $f_+$ is called a maximal function for the kernel in \eqref{4.1aa}. Maximal functions play an important role in the study of Toeplitz kernels, as they determine the kernel uniquely and can be used as test functions for various properties
\cite{CP18}.

We now study the existence of minimal kernels of Toeplitz operators with possibly unbounded symbols.
We start by considering a class of symbols related to Proposition \ref{prop3}.
Let $f_+ \in H_p^+ \setminus \{0\}$ and let $f_+=I_+O_+$ be an inner-outer factorization ($I_+$ inner, $O_+$ outer). Given an outer function $Q$ such that $Q$ is bounded away from zero and $Qf_+ \in H_p^+$, consider the class of symbols $h$ such that $|h|< c|Q|$ for some $c>0$. We denote this class by $\sigma_Q$. We have then $f_+ \in D(T_h)$.

We denote by $K_{\min}^Q(f_+)$ the minimal kernel of Toeplitz operators with symbol in the class $\sigma_Q$. Recall from Section 3 that, if $f_+\in\ker T_g$, then $g$ satisfies  \eqref{2} and we can choose $Q$ according to Proposition \ref{prop3}., with $g\in\sigma_Q$. In the next theorem we show that this minimal kernel exists and we associate it with a symbol in $\sigma_Q$.

%Theorem 4.1
 \begin{thm} \label{thm4.1} For $f_+ \in H_p^+ \setminus \{0\}$,
 \[K_{\min}^Q(f)= \ker T_{\overline{Q}\frac{\overline{I_+}\overline{O_+}}{O_+}}.\]
 \end{thm}

\beginpf Since $Qf_+ \in H_p^+$, we have that $Qf_+=I_+(O_+Q)$ is an inner-outer factorization and

\beq \label{4.3a}
\overline{Q}\frac{\overline{I_+}\overline{O_+}}{O_+}(f_+)=\overline{Q}\frac{\overline{I_+}\overline{O_+}}{O_+}(I_+O_+)=\overline{Q}\overline{O_+} \in H_p^-,
\eeq
because the left-hand side of \eqref{4.3a} is in $L_p$ and the right-hand side is the conjugate of a function in $\Nev_+$ ($O_+^{-1}$ is outer in $H^+_{\infty}$). Therefore
\beq f_+ \in \ker T_{\overline{Q}\frac{\overline{I_+}\overline{O_+}}{O_+}}=\ker T_k \quad (\text{with } k=\overline{Q}\frac{\overline{I_+}\overline{O_+}}{O_+}).\eeq
Now we have to show that any other Toeplitz kernel $\ker T_h$, with $h \in \sigma_Q$ and such that $f_+ \in \ker T_h$, contains $\ker T_k$, i.e.,
$\ker T_k \subset \ker T_h.$
First note that if $h \in \sigma_Q$ and $f_+ \in \ker T_h$, then
$Q^{-1}h \in L_{\infty}$ and $Qf_+ \in \ker T_{Q^{-1}h}.$
So, by \eqref{4.1aa},
\beq \label{4.5a}\ker T_{Q^{-1}h}\supset K_{\min}(Qf_+) = \ker T_{\frac{\overline{Q}}{Q}\frac{\overline{I_+}\overline{O_+}}{O_+}}= \ker T_{Q^{-1}k}. \eeq
Now,
\[\psi_+ \in \ker T_{k} \Leftrightarrow  \psi_+ \in H_p^+, \quad k \psi_+=\psi_- \in H_p^- \Leftrightarrow \]
\beq \label{4.6a}\Leftrightarrow \psi_+ \in H_p^+, \quad kQ^{-1}(Q\psi_+)=\psi_- \in H_p^-.\eeq
So, if $\psi_+ \in \ker T_{k}$ we have, on the one hand,

\[Q\psi_+=\frac{\psi_+}{Q^{-1}} \in \Nev_+\]
and, on the other hand, $Q\psi_+=Q\psi_-/k=\psi_-/(Q^{-1}k) \in L_p$ because $|Q^{-1}k|=1$ a.e.. Therefore,
$Q\psi_+ \in \Nev_+\cap L_p =H_p^+.$
Now, from \eqref{4.6a} and taking \eqref{4.5a} into account, we have
\[Q\psi_+ \in \ker T_{kQ^{-1}}\subset \ker T_{hQ^{-1}},\]
so we conclude that
\[\psi_+ \in \ker T_k \Rightarrow Q\psi_+ \in \ker T_{hQ^{-1}}\Rightarrow
h\psi_+ \in H_p^-\Rightarrow \psi_+ \in \ker T_h. \]
\endpf

A natural question that arises from Theorem \ref{thm4.1} is the relation between $K_{\min}^Q(f_+)$ and $K_{\min}(f_+)$. Since the class $\sigma_Q$ includes all bounded symbols, we have that
% $K_{\min}^Q(f_+) ~\subset~ K_{\min}(f_+)$, more precisely,
\beq \label{K_min_incl}
K_{\min}^Q(f_+)\subset K_{\min}(f_+)\cap \mathcal{D}_Q
\eeq
where
\beq \label{D_Q}
\mathcal{D}_Q= \{\phi_+ \in H_p^+: Q\phi_+ \in L_p\}.
\eeq
To see that the converse of \eqref{K_min_incl} also holds, let $\phi_+ \in K_{\min}(f_+)\cap \mathcal{D}_Q$. Assuming $f_+=I_+O_+$ there exists $\phi_- \in H_p^-$ such that
\[\frac{\overline{I_+}\overline{O_+}}{O_+}\phi_+=\phi_- \, \Rightarrow  \, \overline{Q}\frac{\overline{I_+}\overline{O_+}}{O_+}\phi_+=\overline{Q}\phi_-
\Rightarrow  \, \underbrace{\overline{Q}\frac{\overline{I_+}\overline{O_+}}{O_+}}_{\in L_p}\phi_+=\underbrace{\frac{\phi_-}{(\overline{Q})^{-1}}}_{\in \Nev_+} \in H_p^-.\]
So $\phi_+ \in \ker T_{\overline{Q}\frac{\overline{I_+}\overline{O_+}}{O_+}}$ and we have the following:
%Corollary 4.2
\begin{cor} \label{cor4.2_a}
With the same assumptions as in Theorem \eqref{thm4.1}, for $\mathcal{D}_Q$ defined by \eqref{D_Q},
\[K_{\min}^Q(f_+)=K_{\min}(f_+)\cap \mathcal{D}_Q.\]
\end{cor}

Let us now study the existence of a minimal kernel without restricting to symbols $h \in \sigma_Q$ for some outer $Q$. We denote this minimal kernel for general symbols in $\sigma_p$, if it exists, by $K_{\min}^*(f_+)$.

Note that for any Toeplitz kernel $T_g$ to have a non-zero kernel, we must have
 $gf_+=f_-$ for some $f_{\pm} \in H_p^{\pm} \setminus \{0\}$,
so $g$ is always of the form $f_-/f_+$. If $f_-= \overline{I}\overline{O_+}$ where $I$ is inner and $O_+$ is outer in $H_p^+$, then
\[\ker T_{f_-/f_+} = \ker T_{\overline{I}\overline{O_+}/f_+} \supset \ker T_{\overline{O_+}/f_+},\]
so, when looking for a minimal kernel  containing $f_+$, it is enough to consider symbols of the form $g=\overline{O_+}/f_+$ where $O_+$ is outer in $H_p^+$ or, equivalently, of the form
\beq\label{s}
g={O_-}\overline{I_+}({O_+})^{-1}, \text{ with}\,\, O_\pm \text{ outer in} \,\,H_p^\pm.
\eeq
We say that $O_-$ is outer in $H_p^-$ if  $\overline{O_-}$ is outer in $ H_p^+$.
We start by studying some kernels of Toeplitz operators with symbols of the form \eqref {s}.

%Proposition 4.3
\begin{prop} \label{prop4.2} Let $g={O_-}\overline{I_+}({O_+})^{-1}$ as in \eqref{s}.
Then
\[\ker T_g=\ker T_{{O_-}\overline{I_+}({O_+})^{-1}}={O_+}/{\overline{O_-}}\,\ker T_{\overline{I_+}{O_-}                       /{\overline{O_-}}}\cap H_p^+={O_+}/{\overline{O_-}}\,K_{\min}(I_+\overline{O_-})\cap H_p^+.\]
\end{prop}
\beginpf For $\phi_{\pm}\in H_p^{\pm}$,
\[g\phi_+=\phi_- \Leftrightarrow \frac{O_-}{f_+} \phi_+=\phi_- \Leftrightarrow \frac{O_-}{I_+ O_+}\phi_+=\phi_- \Leftrightarrow   \underbrace{\overline{I_+}\frac{O_-}{\overline{O_-}}}_{\in \mathcal{G} L_{\infty}} \underbrace{\frac{\overline{O_-}}{O_+}\phi_+}_{\in \Nev_+\cap L_p=H_p^+}=\phi_- \Leftrightarrow\]
\[\Leftrightarrow \frac{\overline{O_-}}{O_+}\phi_+ \in \ker T_{\overline{I_+}\frac{O_-}{\overline{O_-}}} \Leftrightarrow \phi_+ \in \frac{O_+}{\overline{O_-}} \ker T_{\overline{I_+}\frac{O_-}{\overline{O_-}}}. \]\endpf

%prop 4.4
\begin{prop} Let $f_+ \in H_p^+ \setminus \{0\}$ and let $f_+=I_+O_+$ be an inner-outer factorization. Then $f_+\in \ker T_{{O_-}\overline{I_+}({O_+})^{-1}}$  if $O_-$ is an outer function in $H_p^-$ and, for  any $h$ such that $f_+ \in \ker T_h$, we have
\[\ker T_h \supset \ker T_{{O_-}\overline{I_+}({O_+})^{-1}}\cap D(T_h).\]
\end{prop}
\beginpf
Suppose that $f_+ \in \ker T_h$, i.e., $hf_+ \in H_p^-$. For any $\phi_+ \in \ker T_{O_-/f_+} \cap D(T_h),$ by Proposition \ref{prop4.2} we have
$\phi_+ =O_+/\overline{O_-}k_+$ with $k_+ \in K_{\min}(I_+\overline{O_-})=\ker T_{\overline{I_+}O_-/\overline{O_-}}$ and
\[\underbrace{h\phi_+}_{\in L_p}=(\underbrace{hI_+O_+}_{\in H_p^-})(\underbrace{\overline{I_+}\frac{O_-}{\overline{O_-}}k_+}_{\in H_p^-})\frac{1}{O_-}\in \overline{\Nev_+}\cap L_p = H_p^-\]
\endpf

If $\ker T_{{O_-}\overline{I_+}({O_+})^{-1}}\subset D(T_h)$ for all $h$ such that $I_+O_+=f_+\in D(T_h)
$, then
\[
\ker T_{{O_-}\overline{I_+}({O_+})^{-1}}=K^*_{\min}(f_+).
\]
As we show next, this holds if $I_+$ is a finite Blaschke product and $O_-$ is a square rigid function in $H_p^-$.

We say that an outer function $\psi_+ \in H_p^+$ is square rigid if and only if $\psi_+^2$ is rigid in $H^+_{p/2}$ (see, e.g., \cite{CP14}).
(A rigid function   $f_+\in H_q^+ \setminus \{0\}$ is one such that every function $g_+ \in H_q^+$ with
$g_+/f_+>0$ a.e.\ on $\RR$ satisfies $g_+=\lambda f_+$ for some $\lambda \in \RR^+$.)
It can be shown that $\span \{\psi_+\}$ is a Toeplitz kernel in $H_p^+$ if and only if $\psi_+$ is square rigid and, in that case,
\[K_{\min}(\psi_+)=\ker T_{\overline{\psi_+}/\psi_+}=\span \{\psi_+\}\]
(\cite{CP14, Sarason}). If $\psi_- \in H_p^-$, we say that $\psi_-$ is square rigid if and only if $\overline{\psi_-}$ is square rigid in $H_p^+$.

%Prop 4.4
\begin{prop} \label{prop4.3}Let $g={O_-}\overline{B}({O_+})^{-1}$ where $B$ is a finite Blaschke product, $O_+ \in H_p^+$ is outer and $O_-$ is a square rigid function in $H_p^-$. Then, if $B$ is a constant,
\beq \label{4.1ab}
\ker T_g=\ker T_{{O_-}({O_+})^{-1}}=\span\{O_+\}\eeq
and, if $B$ is not constant,
\beq \label{4.1a}
\ker T_g=\ker T_{{O_-}\overline{B}({O_+})^{-1}}=\span\{O_+\}\oplus \frac{\xi-z_1}{\xi+i}R_+O_+\span \{r^j, j=0, \cdots, k-1\}\eeq
where $z_1$ is any one of the zeroes of $B$, $k$ is the number of zeroes of $B$ and $R_+$ is rational, invertible in $H^+_{\infty}$, such that $B=R_-r^kR_+$ with
 $R_-=\overline{R_+^{-1}}$.
\end{prop}
\beginpf If $B$ is a constant we may assume it to be $1$.
We have
\[\ker T_{O_-/O_+}={O_+}/{\overline{O_-}}\ker T_{O_-/\overline{O_-}}\cap H_p^+={O_+}/{\overline{O_-}}\span\{\overline{O_-}\}=\span\{O_+\}.\]
 On the other hand, it is easy to see that if $B$ is a finite Blaschke product of degree $k\geq 1$, then we can write $B=R_-r^kR_+$ with $R_{\pm}$ rational and invertible in $H^{\pm}_{\infty}$, such that $R_-=\overline{R_+^{-1}}$ (\cite{5GK,8MP}). We have then $K_B=\overline{R_-^{-1}}K_{r^k}=\frac{R_+}{\xi+i}\span \{r^j, j=0, \cdots,k-1\}$.

% \textcolor[rgb]{0,0,1}{Which paper should we cite here?\\}

Now, from Proposition \ref{prop4.2}, we have that
$\ker T_{O_-\bar{B}({O_+})^{-1}}=\frac{O_+}{\overline{O_-}}\ker T_{\overline{B}O_-/\overline{O_-}}\cap H_p^+$
where
$\ker T_{\overline{B}O_-/\overline{O_-}}=\ker T_{O_-/\overline{O_-}}\oplus \overline{O_-}(\xi-z_1)K_B$
by Theorem 6.7 in \cite{CP14}. Therefore
\begin{eqnarray*}
\ker T_{O_-/f_+} & =& \frac{O_+}{\overline{O_-}}\left(\ker T_{O_-/\overline{O_-}}\oplus \overline{O_-}(\xi-z_1)R_+K_{r^k}\right)\cap H_p^+ \\
&= &\frac{O_+}{\overline{O_-}}\left( \span \{\overline{O_-}\} \oplus \overline{O_-}(\frac{\xi-z_1}{\xi+i})R_+\span \{r^j, j=0, \cdots,k-1\}\right)\cap H_p^+ \\
&= & \span \{O_+\} \oplus (\frac{\xi-z_1}{\xi+i})O_+R_+\span \{r^j, j=0, \cdots,k-1\}.
\end{eqnarray*}
\endpf

%\textcolor[rgb]{0,0,1}{Why is the next result in red?\\}

\begin{cor} If $g=O_-I_+(O_+)^{-1},$ where $O_+ \in H_p^+$ is outer, $O_- \in H_p^-$ is square rigid and $I_+$ is a non-constant inner function, then $\ker T_g=\{0\}$.
\end{cor}
\beginpf
It is easy to see that, as in the case of bounded symbols \cite{CMP},
\[\ker T_{I_+g}\subsetneqq \ker T_g.\]
Since $\ker T_{O_-(O_+)^{-1}}=\span \{O_+\}$, we must have $\ker T_{I_+O_-(O_+)^{-1}}=\{0\}$.
\endpf

From Proposition \ref  {prop4.3} and \eqref{4.1aa} we see that, if $O_-$ is square rigid in $H_p^-,$ then $\ker T_{O_-\bar{B}({O_+})^{-1}}$ does not depend on $O_-$. We have thus the following:

%%%%%Corollary 4.5
\begin{cor} \label{cor4.4} If $g={O_-}\overline{B}({O_+})^{-1})$ where $B$ is a finite Blaschke product, $O_+ \in H_p^+$ is outer and $O_-$ is a square rigid function in $H_p^-$, then  $\ker T_g$ does not depend on $O_-$ and
\begin{eqnarray*}
\ker T_{O_-\bar B({O_+})^{-1})} %&= & \span \{O_+\} \oplus (\frac{\xi-z_1}{\xi+i})O_+R_+\span \{r^j, j=0, \cdots,k-1\} \\
%&= &
=O_+(\xi+i) K_{rB}
\end{eqnarray*}
\end{cor}
\beginpf If $\phi_+ \in \ker T_{O_-/f_+} $ then $\frac{\phi_+}{O_+}\in H^+_{\infty}$ and $\frac{\phi_+}{(\xi+i)O_+}\in H_p^+$.
On the other hand
%\[\frac{O_-}{f_+}\phi_+=
$O_-\frac{\overline{B}\phi_+}{O_+}=\phi_- \in H_p^-,$
so
\[\frac{1}{\xi-i}\frac{\overline{B}\phi_+}{O_+}=\frac{\phi_-}{O_-}\in \,\overline{\Nev_+} \cap L_p=H_p^-.\]

Therefore
\[\overline{r}\overline{B}\frac{\phi_+}{O_+(\xi+i)}=\frac{1}{\xi-i}\overline{B}\frac{\phi_+}{O_+}\in H_p^-,\]
which means that $\frac{\phi_+}{O_+(\xi+i)}\in \ker T_{\overline{r}\overline{B}}=K_{rB}$, and we conclude that \[\ker T_{{O_-}\overline{B}({O_+})^{-1})} \subset O_+(\xi+i)K_{rB}. \]
Conversely, if $\phi_+ \in O_+(\xi+i)K_{rB}$ then $\phi_+ \in H_p^+, \phi_+/O_+ \in H^+_{\infty}$ and
\[\overline{r}\overline{B}\frac{\phi_+}{O_+(\xi+i)}=\phi_- \in H_p^-.\]
Thus
\[\frac{O_-}{f_+}\phi_+=\underbrace{\frac{\overline{B}O_-}{O_+}\phi_+}_{\in L_p}=\underbrace{\frac{O_-}{1/(\xi-i)}}_{\in \overline{\Nev_+}}\in H_p^-,\]
so $\phi_+ \in \ker T_{{O_-}\overline{B}({O_+})^{-1})}$.
\endpf

%Corollary4.5 A

%Theorme 4.6
\begin{thm}\label{thm3} If $f_+=BO_+$ where $B$ is a finite Blaschke product, then
\[K^*_{\min}(f_+)=\ker T_{{O_-}\overline{B}({O_+})^{-1})}= O_+(\xi+i) K_{rB}\]
where $O_-$ is any square rigid outer function in $H_p^-$.
\end{thm}

\beginpf Obviously, $f_+\in \ker T_{O_-/f_+}$. Suppose now that $f_+ \in \ker T_h$; we want to show that $\ker T_h \supset \ker T_{O_-/f_+}$.

Since $f_+=BO_+,$ by near invariance (\cite{CP14}) we also have $O_+ \in \ker T_h$. On the other hand
\[hBO_+=f_- \Leftrightarrow hR_-r^kR_+O_+=f_-\Leftrightarrow\]
\[\forall_{j \in \{0,\cdots, k-1\}} \,,\, h \,\,\underbrace{\left(\frac{\xi-z_1}{\xi+i}R_+r^jO_+\right)}_{\in H_p^+}=\underbrace{\frac{\xi-z_1}{\xi+i}R_-^{-1}r^{-k+j+1}f_-}_{\in H_p^-}.\]

By Proposition \ref{prop4.3}, we conclude that $\ker T_{O_-/f_+} \subset \ker T_h$, so $K^*_{\min}(f_+)=\ker T_{O_-/f_+}$, and the remaining equalities follow from $O_-/f_+\phi_+=\phi_-$ with $\phi_{\pm}\in H_p^{\pm}$.

Then
\beq \label{4.3aa}
\frac{1}{(\xi-i)f_+}\phi_+=\frac{1/(\xi-i)\phi_-}{O_-}
\eeq
where the left hand side belongs to $L_p$ because, by Proposition \ref{prop4.3},
\[\frac{1}{(\xi-i)f_+}\phi_+=\frac{1}{(\xi-i)BO_+}(O_+\widetilde{\phi_+})=B\frac{\widetilde{\phi_+}}{\xi-i}\]
with $\widetilde{\phi_+} \in H^+_{\infty}$, while the right-hand side of \eqref{4.3aa} is in $\overline{\Nev_+}$. Therefore $\frac{1}{(\xi-i)f_+}\in L_p\cap \overline{\Nev_+}=H_p^-$ and we conclude that $\ker T_{O_-/f_+} \subset \ker T_{\frac{1}{(\xi-i)f_+}}$.

The converse can be shown analogously.

Finally, the second equality in this theorem follows from Corollary \ref{cor4.4}.
\endpf

%Corollary 4.7
\begin{cor}
\label{cor:4.9}
A function $f_+ \in H_p^+$ has an inner factor that is a finite Blaschke product of degree $k$ if and only if $\dim (K^*_{\min}(f_+))=k+1$.
\end{cor}

\beginpf By the previous theorem, if $f_+=I_+O_+$ where $I_+$ is a finite Blaschke product of degree $k$, then $K^*_{\min}(f_+)$  has dimension $k+1$. Conversely, if $\dim (K^*_{\min}(f_+))=k+1$, by the property of $\eta$-near invariance of kernels for $\eta \in H_{\infty}^-$ (\cite{CP14}), the inner factor of $f_+$ must be a finite Blaschke product of degree $k$.\endpf

As was mentioned before, $f_+ \in H_p^+$ is square rigid if and only if $K_{\min}(f_+)=\span\{f_+\}$. The next result (proved similarly to Corollary \ref{cor:4.9}) provides an analogous description for outer functions in terms of minimal kernels.

%Corollary 4.7
\begin{cor} A function $f_+ \in H_p^+$ is outer if and only if $K^*_{\min}(f_+)=\span \{f_+\}$.
\end{cor}

%Section 5 - $(p,q) $ - factorization

\section{$(j,s) $ - factorization} \label{sec_pq_fact}
\label{sec:5}

Wiener--Hopf $p$-factorization, presented in Section 2, is of the type considered above. However, even for simple piecewise continuous symbols, it may not exist. In that case the range of the Toeplitz operator is not closed, so the operator is not Fredholm (nor invertible), but the question of describing its kernel and the kernel of its adjoint still stands. In this section we define a more general type of factorization which will allow us to describe the kernels of Toeplitz operators whose symbols may not admit a $p$-factorization, in particular those with piecewise continuous symbols.

%Definition 5.1
\begin{defn}
%Let $g, g^{-1} \in \BB_s$, $s \geq 1$.
A representation of the form
\beq\label{eq_fact}
g=g_-r^kg_+, \quad \text{with } k \in \ZZ,
\eeq
\beq
g_- \in \HH_j^-, \qquad g_-^{-1} \in \HH_s^-,
\eeq
\beq
g_+ \in \HH_s^+, \qquad g_+^{-1} \in \HH_j^+,
\eeq
where $j,s \in ]1, +\infty[\,,\,1/j+1/s\leq 1$, is called a \emph{$(j,s)$-factorization} of $g$ with index $k$.
\end{defn}
If $k=0$, the factorization is said to be \emph{canonical}. If $j=p\,,\,s=p'$ then \eqref{eq_fact} is a {$p$-factorization} of $g$.%;  if $g_+ \notin \HH_{j'}^+$ or $g_-^{-1} \notin \HH_{j'}^-$, we say that \eqref{eq_fact} is a \emph{strict $(j,s)$-factorization}.

It is easy to see that if $g$ admits a $p$-factorization, then it is unique up to non-zero multiplicative constants in $g_{\pm}$ \cite{5GK, 8MP}. However, this is not true for general $(j,s)$ factorizations. Assume that $g$ has two $(j,s)$ - factorizations
\beq \label{two_fact}
g=g_-r^kg_+, \qquad g=\tilde{g}_-r^{\tilde{k}}\tilde{g}_+,
\eeq
where $k\geq\tilde k$ %$g_-^{-1} \notin \HH_{j'}^-$ and $g_+ \notin \HH_{j'}^+$.
and $s\neq j'$. Then
\beq \label{g+g+}
g_+\tilde{g}_+^{-1}r^{k-\tilde{k}}=\tilde{g}_-g_-^{-1}.
\eeq
If $k>\tilde{k}$, the left-hand side of \eqref{g+g+} belongs to $\BB_l^+$, with $\frac{1}{l}=\frac{1}{j}+\frac{1}{s}<1$ , and vanishes at $i$; the right-hand side belongs to $\BB_l^-$, so both sides must be equal to a polynomial $A(\xi-i)$, $A \in \CC$. So, if $k-\tilde{k}\geq 2$, we must have $A=0$, which is impossible.

If $k \leq \tilde{k}$, then from \eqref{g+g+} we get that
\beq
g_+^{-1}\tilde{g}_+r^{\tilde{k}-k}=\tilde{g}_-^{-1}g_-=\tilde{A}\xi+\tilde{B}, \qquad \tilde{A}, \tilde{B} \in \CC,
\eeq
and if $\tilde{k}-k\geq 2$ we must have $\tilde{A}= \tilde{B}=0$, which is impossible. So we must have $0 \leq k-\tilde{k}\leq 1$, and we have the following:

%Proposition 5.2

\begin{prop} \label{prop3.2}If $g$ admits two $(j,s)$-factorizations as in \eqref{two_fact}, then
\begin{itemize}
  \item[(i)] if $g_+ \in \HH_{j'}^+$, $g_-^{-1} \in \HH_{j'}^-$, we have $k=\tilde{k}$ and $g_+\tilde{g}_+^{-1}=\tilde{g}_-g_-^{-1}=C \in \CC \setminus \{0\}$;
  \item[(ii)] if $g_+ \in \HH_{s}^+$ or $g_-^{-1} \in \HH_{s}^-$ with $s\neq j'$ then $0\leq k-\tilde{k}\leq 1$ and
\beq \label{3.8}
\text {if} \quad k=\tilde{k}, \,\, \text{then}\,\, g_+\tilde{g}_+^{-1}=\tilde{g}_-g_-^{-1}=C; \eeq
\beq \label{3.9}
\text {if} \quad k-\tilde{k}=1, \,\, \text{then}\,\,  \tilde{g}_+=\frac{1}{C}\frac{g_+}{\xi+i}\,\,, \,\, \tilde{g}_-=C(\xi-i)g_- \,,\,\, \text {with} \,\,C \in \mathbb C \setminus \{0\}.\eeq
    %\item if $k-\tilde{k}=-1$, \beq \label{3.10} \tilde{g}_+=C(\xi+i)g_+, \qquad \tilde{g}_-=\frac{1}{C}\frac{g_-}{\xi-i}, \eeq

\end{itemize}
\end{prop}

\beginpf
Only the second part of (ii) is left to prove. Assume that $g_+ \in \HH_{s}^+$, $\tilde{g}_- \in \HH_{j}^-$; then
\[g_+\tilde{g}_+^{-1}r^{k-\tilde{k}}=g_-^{-1}\tilde{g}_-.\]
If $k=\tilde{k}$, we have $g_+\tilde{g}_+^{-1}=g_-^{-1}\tilde{g}_-=A \xi+B$, with $A, B \in \mathbb{C}$, by \eqref{BB_inters} and, analogously,
\[g_+^{-1}\tilde{g}_+=g_-\tilde{g}_-^{-1}=\tilde{A} \xi+\tilde{B},\]
with $\tilde{A}, \tilde{B} \in \mathbb{C}$. So we must have $A=\tilde{A}=0$, $B=\frac{1}{\tilde{B}}\neq 0$.

If $k-\tilde{k}=1$, then $g_+\tilde{g}_+^{-1}r=g_-^{-1}\tilde{g}_-=A(\xi-i)$, and it follows that \eqref{3.9} holds.
%If $k-\tilde{k}=-1$, we have $g_+^{-1}\tilde{g}_+r=g_-\tilde{g}_-^{-1}=A( \xi-i),$ and \eqref{3.10} holds.
\endpf

%Corollary 5.3
\begin{cor} If $g=g_-g_+$ and $g=\tilde{g}_-\tilde{g}_+$ are two canonical $(j,s)$-factorizations, then the factors are unique up to a non-zero multiplicative constant.
\end{cor}
Note that the estimate $k-\tilde{k}\leq 1$ in Proposition \ref{prop3.2} (ii), is optimal. For example, $\left(\frac{\xi-i}{\xi+i}\right)^{1/2}_{\infty}$ admits two $(j,s)$ - factorizations, for $j,s>2$, with $k-\tilde k=1$:
\[\left(\frac{\xi-i}{\xi+i}\right)^{1/2}_{\infty}=\sqrt{\xi-i}\frac{1}{\sqrt{\xi+i}}=\frac{1}{\sqrt{\xi-i}}\frac{\xi-i}{\xi+i}\sqrt{\xi+i}\]
%(with appropriate branches).

%
In the case where $g$ admits two factorizations,
 with respect to different pairs of function spaces, we have the following:

%Proposition 5.4
\begin{prop} \label{prop3.4} Let $g=g_-r^kg_+$ be a $(p,q)$-factorization and $g=\tilde{g}_- r^{\tilde{k}} \tilde{g}_+$ be a $(j,s)$~-~factorization. If $j=q'$, then $k\leq \tilde{k}$; if $s=p'$, then $k\geq \tilde{k}$.
\end{prop}

\beginpf Suppose that $j=q'$. From $\tilde{g}_+^{-1}g_+r^{k-\tilde{k}}=g_-^{-1}\tilde{g}_-$, we see that if $k>\tilde{k}$ then both sides of the previous equality must be zero, which is impossible, because the left-hand side belongs to $\BB_1^+$ and vanishes at $i$, while the right-hand side belongs to $\BB_1^-$. So we have $k\leq \tilde{k}$. The second part is proved analogously.
\endpf
For example, we have a $(p,2)$ - factorization, with $p>2$, $\left(\frac{\xi-i}{\xi+i}\right)^{1/2}_{\infty}=\sqrt{\xi-i}\frac{1}{\sqrt{\xi+i}}$, with index 0, and a $(2,s)$ - factorization, with $s>2$, $\left(\frac{\xi-i}{\xi+i}\right)^{1/2}_{\infty}=\frac{1}{\sqrt{\xi-i}}\frac{\xi-i}{\xi+i}\sqrt{\xi+i}$, with index $1>0$.

%Corollary 5.5
\begin{cor}\label{cor3.5} If $g$ admits a $(j,s)$-factorization $g=g_-r^kg_+$, and a $(s',j')$ - factorization $g=\tilde{g}_-r^{\tilde{k}}\tilde{g}_+$, then $k=\tilde{k}$ and we have $g_-\tilde{g}_-^{-1}=g_+^{-1}\tilde{g}_+=C \in \mathbb{C} \setminus \{0\}$.
\end{cor}

\beginpf From Proposition \ref{prop3.4}, it follows that $k=\tilde{k}$ and thus $g_-g_+=\tilde{g}_-\tilde{g}_+$. Since $g_-\tilde{g}_-^{-1}=g_+^{-1}\tilde{g}_+$ where the left-hand side is in $\BB_1^-$ and the right-hand side is in $\BB_1^+$, both sides are equal to a non-zero constant.
\endpf

We will also need the following.

%Proposition 5.6
\begin{prop} \label{prop3.6} If $g$ admits a canonical $(j, p')$-factorization $g=\tilde{g}_-\tilde{g}_+$ where $\tilde{g}_- \notin \HH_p^-$ or $\tilde{g}_+^{-1} \notin \HH_p^+$, and $g$ also admits a $p$-factorization $g=g_-r^kg_+$, then $k=1$ and we have $g_-=Br_-\tilde{g}_- \in H_j^-$, $g_+^{-1}=Br_+\tilde{g}_+^{-1} \in H_j^+$, $B \in \mathbb{C} \setminus \{0\}$.
\end{prop}
\beginpf From $g=\tilde{g}_-\tilde{g}_+=g_-r^kg_+$ and Proposition \ref{prop3.4}, it follows that $k\geq 0$. If $k=0$,  we get $g_-\tilde{g}_-^{-1}=\tilde{g}_+g_+^{-1}$ where $g_-\tilde{g}_-^{-1} \in \BB_1^-$ and $\tilde{g}_+g_+^{-1} \in \BB_1^+$. From \eqref{BB_inters}, we get $g_-\tilde{g}_-^{-1}=\tilde{g}_+g_+^{-1}=C \neq 0$  implying that $\tilde{g}_-\in \HH_p^-$ and $\tilde{g}_+^{-1}\in \HH_p^+$, which is impossible. So we must have $k>0$. Suppose that $k>1$; then we have
\beq \label{3.11} r^kg_-\tilde{g}_-^{-1}=g_+^{-1}\tilde{g}_+=\frac{P_k}{(\xi+i)^k}
\eeq
where $P_k$ is a polynomial of degree smaller or equal to $k$. From \eqref{3.11} we see that $P_k$ cannot have zeroes in $\mathbb{C}^+$ nor in $\mathbb{C}^-$, and from
\[\tilde{g}_+^{-1}=\frac{(\xi+i)^k}{P_k}g_+^{-1}, \quad k>1\]
we see that $P_k$ cannot have zeroes in $\mathbb{R}$ either; thus $P_k=C \in \mathbb{C}\setminus \{0\}$. In that case, we would have
\[g_+\frac{\tilde{g}_+^{-1}}{\xi+i}=\frac{1}{C}(\xi+i)^{k-1}, \quad k-1>0\]
where the left-hand side belongs to $\HH_s^+,$ $\frac{1}{s}=\frac{1}{p'}+\frac{1}{j}<1$, which is impossible for the right-hand side if $k-1>0$.

Therefore we can only have $k=1$ and in that case
\[g_-rg_+=\tilde{g}_-\tilde{g}_+ \Rightarrow rg_-\tilde{g}_-^{-1}=\tilde{g}_+g_+^{-1}=\frac{A\xi+B}{\xi+i}.\]
Thus $A\xi+B$ cannot have zeroes in $\mathbb{C}^{\pm}$; on the other hand we have $\tilde{g}_+^{-1}g_+=\frac{\xi+i}{A\xi+B}$ where the left-hand side belongs to $\BB_s^+$, $s>1$, so $A\xi+B$ cannot have zeroes on $\mathbb{R}$ either.
Therefore $\tilde{g}_+g_+^{-1}=\frac{B}{\xi+i}$ and $\tilde{g}_-^{-1}g_-=\frac{B}{\xi-i}$.
\endpf

%Corollary 5.7
\begin{cor} \label{cor3.7} If $g$ admits a canonical $(j, p')$-factorization $g=\tilde{g}_-\tilde{g}_+$ where $\tilde{g}_- \notin \HH_p^-$ or $\tilde{g}_+^{-1} \notin \HH_p^+$, then $g$ does not admit a canonical $p$-factorization.
\end{cor}

%Example 5.8
\begin{ex} Let $g=\left(\frac{\xi-i}{\xi+i}\right)^{1/2}_{\infty}$, which has a canonical $(p,2)$-factorization $(p>2)$
\[\left(\frac{\xi-i}{\xi+i}\right)^{1/2}_{\infty}=\sqrt{\xi-i}\frac{1}{\sqrt{\xi+i}}\]
(with appropriate branches). Since $g$ is 2-singular, it does not admit a Wiener--Hopf 2-factorization; by Corollary \ref{cor3.7}, it also does not admit a canonical 2-factorization. %(in which condition (iii) of generalized factorization is not considered)
Moreover, if $g$ admitted a (non-canonical) 2-factorization,
by Proposition \ref{prop3.6} its index would be 1; so we would have
\[g_-rg_+=\sqrt{\xi-i}\frac{1}{\sqrt{\xi+i}}, \quad g_-^{\pm 1} \in \HH_2^-, \, g_+^{\pm 1} \in \HH_2^+.\]
This implies that
\[rg_+\sqrt{\xi+i}=g_-^{-1}\sqrt{\xi-i}=\frac{A}{\xi+i}, \quad A\in \mathbb{C},\]
because $g_+\sqrt{\xi+i} \in \HH_p^+$ for any $p>2$, and $g_-^{-1}\sqrt{\xi-i} \in \HH_p^-$. Therefore $Ag_-=(\xi+i)\sqrt{\xi-i}$, and we have $Ag_- \in H_j^- \Leftrightarrow A=0.$ We conclude that $\left(\frac{\xi-i}{\xi+i}\right)^{1/2}_{\infty}$ does not admit any 2-factorization.
\end{ex}

%Section 6 - Properties of $(j.s)$-factors
\section{Properties of $(j,s)$-factors}
\label{sec:6}

We start by considering the case where $j=p\,,\,s=p'$.

%THeorem 6.1
\begin{thm} Let $1<p<\infty$. If $g_{\pm}$ is such that
\beq \label{4.1} g_{\pm} \in \HH_{p'}^{\pm}, \quad g_{\pm}^{-1} \in \HH_{p}^{\pm},
\eeq
then $\log g_{\pm} \in \HH_2^{\pm}$.
\end{thm}

\beginpf Let $g_+$ satisfy \eqref{4.1}. Then its analytic extension to the upper half plane $\mathbb{C}^+$ is such that $g_+^{\pm 1}(z) \neq 0,$ for all $z \in \mathbb{C}^+$. Thus we can define an analytic branch of $\log g_+$ in $\mathbb{C}^+$ (\cite{1A}). On the other hand for
\beq \label{4.2} \tilde{g}_+=r_+^{1-2/p'}g_+
\eeq
we have
\[r_+^{2/p'}\tilde{g}_+=r_+g_+ \in H_{p'}^+, \quad r_+^{2/p}\tilde{g}_+^{-1}=r_+g_+^{-1} \in H_{p}^+.\]

Therefore, defining $\widetilde{G}_+(w)=\tilde{g}_+\left(i\frac{1+w}{1-w}\right)$, $w \in \mathbb{D}$, we have $\widetilde{G}_+ \in H^{p'}(\mathbb{D})$, $\widetilde{G}_+^{-1} \in H^{p}(\mathbb{D})$ (\cite{3Dur, 7Koosis}). Let $\log \widetilde{G}_+$ be analytic in $\mathbb{D}$; we have $\log \widetilde{G}_+ \in H^2(\mathbb{D})$ if
\beq \label{4.3} \sup_{0<r<1}\int_{-\pi}^{\pi}|\re\, \log \widetilde{G}_+(re^{i\theta})|^2\,d\theta < \infty.
\eeq
Defining, for each $r \in ]0,1[$,
\[I_1=\{\theta \in [-\pi,\pi]: |\widetilde{G}_+(re^{i\theta})|\geq 1\}\]
\[I_2=\{\theta \in [-\pi,\pi]: 0<|\widetilde{G}_+(re^{i\theta})|< 1\},\]
we have that
\[\int_{-\pi}^{\pi}|\re\, \log \widetilde{G}_+(re^{i\theta})|^2\,d\theta= \int_{-\pi}^{\pi}\log^2 |\widetilde{G}_+(re^{i\theta})|\,d\theta=\]
\[=\int_{I_1}\frac{\log^2 |\widetilde{G}_+(re^{i\theta})|}{|\widetilde{G}_+(re^{i\theta})|^{p'}}|\widetilde{G}_+(re^{i\theta})|^{p'}\,d\theta+\int_{I_2}\frac{\log^2 |\widetilde{G}_+^{-1}(re^{i\theta})|}{|\widetilde{G}_+^{-1}(re^{i\theta})|^{p}}|\widetilde{G}_+^{-1}(re^{i\theta})|^{p}\,d\theta\]
\[\leq M_1\int_{I_1}|\widetilde{G}_+(re^{i\theta})|^p\, d\theta+M_2\int_{I_2}|\widetilde{G}_+^{-1}(re^{i\theta})|^p\,d\theta\]
\[\leq M_1\int_{-\pi}^{\pi}|\widetilde{G}_+(re^{i\theta})|^p\, d\theta+M_2\int_{-\pi}^{\pi}|\widetilde{G}_+^{-1}(re^{i\theta})|^p\,d\theta,\]
where $M_1$ and $M_2$ are positive constants and we used the fact that $\frac{\log^2 x}{x^p}$ is a bounded function in $]1,+\infty[$, for $p>1$.

Since $\widetilde{G}_+ \in H^{p'}(\mathbb{D})$, $\widetilde{G}_+^{-1} \in H^p(\mathbb{D})$, it follows that \eqref{4.3} holds, and we conclude that $\log \widetilde{G}_+ \in H^2(\mathbb{D})$. Thus $\log g_+ \in \HH_2^+$ and, since $\log r_+^{1-2/p'} \in \HH_2$ (with appropriate branches), we have from \eqref{4.2} that $\log g_+ \in \HH_2^+$.

The result regarding $g_-$ is proved analogously.
\endpf

%Corollary 6.2
\begin{cor} \label{cor4.2} If  $g=g_-g_+$ where $g_{\pm}$ satisfy \eqref{4.1} for some $p>1$, then $\log g \in \mathcal{L}_2.$
\end{cor}
Let now $\Pi^{\pm}$ be the bounded complementary projections in $\mathcal{L}_2$ defined by
\[\Pi^{\pm}:\mathcal{L}_2 \rightarrow \HH_2^{\pm}, \qquad \Pi^{\pm}\varphi=r_+^{-1} P^{\pm}(r_+\varphi).\]
Note that $\Pi^- \varphi$ can be expressed equivalently as
\[ \Pi^{-}\varphi=r_-^{-1}P^-r_-\varphi - \frac{1}{\pi}\int_{\mathbb{R}}\frac{\varphi(t)}{1+t^2}\,dt, \quad \varphi \in \mathcal{L}_2. \]

A simple consequence of the previous results is the following.
%Corollary 6.3
\begin{cor} \label{cor4.3} If  $g$ admits a canonical $p$-factorization $g=g_-g_+$ on $\mathbb{R}$, then the factors $g_{\pm}$ are given (up to a multiplicative constant) by
\[g_{\pm}=\exp(\Pi^{\pm}\log g)\]
\end{cor}
\beginpf
If $g$ admits a canonical $p$-factorization $g=g_-g_+$, then by Theorem \ref{4.1} we have $\log g_{\pm} \in \HH_2^{\pm}$ and
$\log g_-+ \log g_+=\Pi^+ \log g + \Pi^-\log g,$
which is equivalent to
\[\underbrace{\log g_--\Pi^-\log g}_{\in \HH_2^-}=\underbrace{\Pi^+\log g - \log g_+}_{\in \HH_2^+}\]
so both sides are equal to a constant $C \neq 0$, and it follows that $g_{\pm}=C \exp(\Pi^{\pm}\log g)$.
\endpf
Obviously, if $g$ admits a $p$-factorization with index $k$, then $g=g_-r^kg_+$ with $g_{\pm}= \exp (\Pi^{\pm} \log g_0)$, $g_0=gr^{-k}$.

Corollary \ref{cor4.3} generalizes a similar result obtained in \cite{14W} for piecewise continuous functions. It might also be obtained by relating a $p$-factorization with respect to $\mathbb{R}$ with a generalized factorization of $g\left(i\frac{1+w}{1-w}\right)$, $w \in \Gamma_0$, in a weighted $L_p$ space of the unit circle $L^p(\Gamma_0,|1-w|^{1-2/p})$ (\cite{CDR, 5GK}), and generalizing to these weighted $L_p$ spaces the existing formulas for the factors in the case of
$L^p(\Gamma)$, where $\Gamma$ is a closed rectifiable contour \cite{5GK, LS}). This has not been done in the published literature, at least explicitly, to the authors' knowledge, and we take here a different approach, by studying the behavior of the functions satisfying condition \eqref{4.1} on $\mathbb{R}$.
%
%\vspace{2mm}

%It is easy to see also the following.

%Proposition 6.4
%\begin{prop}\label{prop4.4} If $b$ has an $(m,l)$ - factorization ($\frac{1}{l}+\frac{1}{m} \leq 1$) $b=b_-r^{k_1}b_+$, and $c$ has a $(j,s) - $ factorization ($\frac{1}{j}+\frac{1}{s} \leq 1$) $c=c_-r^{k_2}c_+$, with $\frac{1}{\tilde{p}}+\frac{1}{\tilde{q}} \leq 1$ where $\frac{1}{\tilde{p}}=\frac{1}{m}+\frac{1}{j}$ and $\frac{1}{\tilde{q}}=\frac{1}{l}+\frac{1}{s}$, then $a=bc$ admits a $(\tilde{p},\tilde{q})$ - factorization $a=(b_-c_-)r^{k_1+k_2}(b_+c_+)$.
%\end{prop}
It may happen that a given function $a$ is known to have a factorization with respect to certain Hardy spaces
%, for example a $p$-factorization (see section 2)
 and, on the other hand, $a$ can be written as the product of two other functions $b$, $c$, each having a particular factorization. The question then arises of how the factorization of $a$ is related to the factorizations of the factors $b$ and $c$.  We have the following:

%Theorem 6.4
\begin{thm}\label{thm4.5} Let $a$ admit a canonical $(p,q)$-factorization $a=a_-a_+$, and let $a=bc$, where $b$ admits an $(m, l)$-factorization $b=b_-b_+$ and $c=c_-c_+$ is a $(j, s)$-factorization, with $\frac{1}{l}+\frac{1}{s} = \frac{1}{p'}$ and $\frac{1}{m}+\frac{1}{j}=\frac{1}{q'}$. Then $a=(b_-c_-)(b_+c_+)$ is a canonical $(p,q)$-factorization of $a$.
\end{thm}
\beginpf We start by remarking that if $\frac{1}{p}+\frac{1}{q}\leq 1$, then $\frac{1}{p'}+\frac{1}{q'}\leq 1$.
We have $a=(b_-c_-)(b_+c_+)$, which is a canonical $(q',p')$ - factorization by Proposition \ref{prop3.4}, and on the other hand $a=a_-a_+$ is a $(p,q)$ - factorization.
From Corollary \ref{cor3.5} it follows that $a_{\pm}=b_{\pm}c_{\pm}$, up to a constant factor.
\endpf
%An immediate consequence of this result is that if  $g$ admits an $(\frac{n}{q'},\frac{n}{p'} )$ - factorization, \textcolor[rgb]{1.00,0.00,0.00}{($n \in \mathbb{N}$ such that $\frac{q'+p'}{n} \leq 1$)}, $g=g_-r^kg_+$, with
%\beq \label{4.4} g_- \in \HH_{n/q'}^-, \,  g_-^{-1} \in \HH_{n/p'}^-, \,  g_+ \in \HH_{n/p'}^+,\,  g_+^{-1} \in \HH_{n/q'}^+,
%\eeq
%and $g^n$ admits a $(p,q)$ - factorization, then $g^n=g_-^nr^{kn}g_+^n$ is a $(p,q)$ - factorization of $g^n$ and we have that
%\beq \label{4.5} g_-^n \in \HH_{p}^-, \,  g_-^{-n} \in \HH_{q}^-, \,  g_+^n \in \HH_{q}^+, \,  g_+^{-n} \in \HH_{p}^+.
%\eeq
%In particular, from Theorem (...$g \in \mathcal{G}C_{\mu}(\mathbb{R})$...) into account we have thus (\cite{5GK})

%\textcolor[rgb]{0,0,1}{Reference?\\}

%\begin{cor} \label{cor4.6} If  $g\in \mathcal{G}C_{\mu}(\dot{\mathbb{R}})$ and $g=g_-r^kg_+$ is a generalized factorization of $g$, then
%\[g^n=g_-^nr^{kn}g_+^n
%\] is a generalized factorization of $g^n$.
%\end{cor}
We have moreover the following.

%Prop 6.5
\begin{prop} \label{cor4.7} If  $g_{\pm}^n \in \HH_p^{\pm}$ for every $n \in \mathbb{N}$, $p>1$, with $g^{\pm}(z)\neq 0$ for all $z \in \mathbb{C}^{\pm}$, then $g_{\pm}^{\alpha} \in \HH_p^{\pm}$, for all $p>1$, $\alpha \in \mathbb{R}$.
\end{prop}
\beginpf Let $n=[\alpha]$ be the integer part of $\alpha$; we have
\beq \label{4.6} g_-^n, \, g_-^{n+1} \in \HH_p^- \text{ for every } p>1.
\eeq
%Let $g_-^{\alpha}$ be analytic in $\mathbb{C}^-$ (note that $g_-$ has a non-vanishing extension to $\mathbb{C}^-$), and define,
Defining, for each $y<0$,
\[\mathbb{R}_1=\left\{x \in \mathbb{R}: 0<|g_-(x+iy)|\leq 1\right\}\]
\[\mathbb{R}_2=\left\{x \in \mathbb{R}: |g_-(x+iy)|> 1\right\},\]
 we have, for all  $y<0$,
\[\int_{\mathbb{R}}\left|\frac{g_-^{\alpha}(x+iy)}{x+iy-i}\right|^p\,dx \leq \int_{\mathbb{R}_1}\left|\frac{g_-^{n}(x+iy)}{x+iy-i}\right|^p\,dx+\int_{\mathbb{R}_2}\left|\frac{g_-^{n+1}(x+iy)}{x+iy-i}\right|^p\,dx\]
\[\leq \int_{\mathbb{R}}\left|\frac{g_-^{n}(x+iy)}{x+iy-i}\right|^p\,dx+\int_{\mathbb{R}}\left|\frac{g_-^{n+1}(x+iy)}{x+iy-i}\right|^p\,dx\]
and from \eqref{4.6} we conclude that $r_-g_-^{\alpha}\in H_p^-$. The result for $g_+^{\alpha}$ can be proved analogously.
\endpf
\begin{cor} \label{cor4.8} With the same assumptions as in Corollary \ref{cor4.7} and $\alpha \in \mathbb{R}$, we have that
\beq \label{4.7} r_{\pm}^{\alpha}g_{\pm} \in H_p^{\pm} \text{ for every } \alpha>\frac{1}{p}.
\eeq
\end{cor}
The assumptions of Theorem \ref{thm4.5} and Corollaries \ref{cor4.7} and \ref{cor4.8} are satisfied, in particular, if $g  \in \mathcal{G}C_{\mu}(\dot{\mathbb{R}})$ and these results provide a partial description of the behavior of the factors $g_{\pm}$ in a neighborhood of any point in $\mathbb{R}$ and at $\infty$. In particular in a neighborhood of $\infty$ they show that although $g_-^{\pm}$, $g_+^{\pm}$ may be unbounded, they ``grow less than'' any positive power of $\xi$.
%We refer to functions satisfying \eqref{4.7} for all $p>1$ as \emph{almost bounded}.
Corollary \ref{cor4.8} also means that  the domain of the Toeplitz operator $T_{g_+}$ (in $H_p^+$) contains all functions of the form $\frac{1}{(\xi+i)^{\alpha}}$ with $\alpha>\frac{1}{p}$, and analogously for  $T_{g_-}$.

%Section 7 - Kernels of Toeplitz operators and $(q,p')$ - factorization
\section{Kernels of Toeplitz operators in $H_p^+$ and $(q,p')$-factorization}
\label{sec:7}

In what follows we assume that $g$, not necessarily bounded, admits a $(q,p')$-factorization ($1<p \leq q <\infty$)
\beq \label{5.1} g=g_-r^kg_+ \text{ with } k \in \mathbb{Z},
\eeq
\beq \label{5.2}g_- \in \HH_q^-, \, \, g_-^{-1} \in \HH_{p'}^-, \, \, g_+ \in \HH_{p'}^+, \, \, g_+^{-1} \in \HH_{q}^+.
\eeq
Note that $1<p \leq q <\infty$ implies that $\frac{1}{q}+\frac{1}{p'}\leq 1.$

%Also recall that $\spam \left\{\frac{1}{\xi+i},\frac{\xi-i}{\xi+i}\frac{1}{\xi+i}, \cdots , \left(\frac{\xi-i}{\xi+i}\right)^{n-1}\frac{1}{\xi+i} \right\}$ is the model space $K_{r^n}$, $n \in \mathbb{N}$.

%Theorem 7.1
\begin{thm} \label{thm5.1} Let $g$ admit a representation \eqref{5.1}-\eqref{5.2}.
\begin{description}
  \item[(i)] If $g_- \in \HH_p^-$, $g_+^{-1}\in \HH_p^+$, then
 $\dim \ker T_g=0$ if $k\geq 0$,
  $\dim \ker T_g=|k|$ if $k< 0$ and, in this case, $\ker T_g=g_+^{-1}K_{r^{|k|}}$.

  \item[(ii)] If $g_- \notin \HH_p^-$ or $g_+^{-1}\notin \HH_p^+$, then
 $\dim \ker T_g=0$ if $k\geq -1$,
 $\dim \ker T_g=|k|-1$ if $k\leq -2$ and, in this case, $\ker T_g=\frac{g_+^{-1}}{\xi+i}K_{r^{|k|-1}}$.
\end{description}
\end{thm}

Before proving Theorem \ref{thm5.1}, note that, by Proposition  \ref{prop3.6}, if the assumptions of (ii) are satisfied and $g$ also admits a $p$-factorization, then the latter has index $k+1.$
\beginpf We have that $\varphi_+ \in \ker T_g$ if and only $\varphi_+ \in \mathcal{D}(T_g)$ and $P^+g\varphi_+=0$, i.e. $\varphi_+ \in H_p^+$ and $g\varphi_+=\varphi_- \in H_p^-.$
Now, for $\varphi_{\pm} \in H_p^{\pm}$,
\beq \label{5.3}  g\varphi_+=\varphi_- \Leftrightarrow r^k \underbrace{g_+\varphi_+}_{\in \HH_1^+}=\underbrace{g_-^{-1}\varphi_-}_{\in \HH_1^-}.
\eeq
If $k\geq 0$, both sides of the last equality must be equal to 0. For $k=-1$, both sides must be equal to a function of the form $\frac{A}{\xi-i}$, with $A \in \mathbb{C}$,
 and it follows that
 $\varphi_-=\frac{A}{\xi-i}g_- \in H_p^-, \qquad \varphi_+=\frac{A}{\xi+i}g_+^{-1} \in H_p^+.$
 Therefore, in the case (i), we have $\ker T_g =g_+^{-1}\spam \left\{\frac{1}{\xi+i}\right\}$ and, in case (ii), we have $A=0$ and $\ker T_g=\{0\}$.
 For $k\leq -2$, both sides of the second equality in \eqref{5.3} must be of the form $\frac{P_{|k|-1}}{(\xi-i)^{|k|}}$, where $P_{|k|-1}$ is a polynomial of degree smaller or equal
 to $|k|-1$. It follows that
 \beq \label{5.4} \varphi_-=\frac{g_-}{\xi-i}\frac{P_{|k|-1}}{(\xi-i)^{|k|-1}}  \in H_p^-, \quad \varphi_+=\frac{g_+^{-1}}{\xi+i}\frac{P_{|k|-1}}{(\xi+i)^{|k|-1}}  \in H_p^+.
 \eeq
 If the degree of $P_{|k|-1}$ is equal to $|k|-1$, then $\frac{(\xi-i)^{|k|-1}}{P_{|k|-1}}$ is analytic and bounded in a neighborhood of $\infty$;
  since $\varphi_- \in H_p^-$  and, from \eqref{5.4}, $\frac{(\xi-i)^{|k|-1}}{P_{|k|-1}}\varphi_-=\frac{g_-}{\xi-i} \in H_q^-$, with $q \geq p$,
  we must have $\frac{(\xi-i)^{|k|-1}}{P_{|k|-1}}\varphi_- \in H_p^-$, and therefore $\frac{g_-}{\xi-i} \in H_p^-$. We conclude analogously that
$\frac{g_+^{-1}}{\xi+i} \in H_p^+$. Thus, in the case (ii) the degree of $P_{|k|-1}$ must be smaller or equal to $|k|-2$.

Conversely, in the case (i) we see that every $\varphi_+$ of the form given in \eqref{5.4} belongs to $\ker T_g$ and $\ker T_g=g_+^{-1}K_{r^{|k|}}$; in the case (ii) we see that every $\varphi_+=g_+^{-1}\frac{P_{|k|-2}}{(\xi+i)^{|k|}}$ belongs to $\ker T_g$ and $\ker T_g=\frac{g_+^{-1}}{\xi+i}K_{r^{|k|-1}}$.
\endpf

It follows from Theorem \ref{thm5.1} that Toeplitz operators in $H_p^+$ with unbounded symbols possessing a $(q, p')$ - factorization, have the same kernel as a Toeplitz operator with bounded symbol; we have $\ker T_g=\ker T_{g_1}$ with $g_1=\overline{g_+^{-1}}r^kg_+$.

\vspace{3mm}

 If $g \in L^{\infty}$, then $T_g$  is bounded in $H_p^+$ and $T^*_g=T_{\bar{g}}$ is bounded in $H_{p'}^+$. In that case the factorization \eqref{5.1}-\eqref{5.2} also allows us to describe $\ker T_g^*$.

%Theorem 7.2
 \begin{thm} \label{thm5.2} Let $g \in L^{\infty}$, such that  \eqref{5.1}-\eqref{5.2} hold.
 \begin{description}
   \item[(i)] If $g_- \in \HH_p^-$, $g_+^{-1}\in \HH_p^+$, we have that
$\dim \ker T^*_g=0$ if $k\leq  0$;
 $\,\dim \ker T^*_g=k$ if $k> 0$ and, in this case, $\ker T^*_g=\overline{g_-^{-1}}K_{r^{k}}$.

   \item[(ii)] If $g_- \notin \HH_p^-$ or $g_+^{-1}\notin \HH_p^+$, we have that
 $\dim \ker T^*_g=0$ if $k\leq 0$;
$\,\dim \ker T^*_g=k$ if $k>0$ and either $g_-^{-1} \notin \HH_{p'}^-$ or $g_+\notin \HH_{p'}^+$;
$\,\dim \ker T^*_g=k+1$ if $k>0$, $g_-^{-1} \in \HH_{p'}^-$ and $g_+\in \HH_{p'}^+$.

 \end{description}
 \end{thm}

\beginpf
[(i)] The result follows from Theorem \ref{thm5.1}, since the existence of a $p$-factorization  for $g$ is equivalent to the existence of a $p'$-factorization for $\bar{g}$, with symmetric indices.

[(ii)] If $\psi_{\pm} \in H_{p'}^{\pm}$ and $\bar{g}\psi_+=\psi_-$, then we have
           \beq \label{5.5} r^{-k} \bar{g}_-\psi_+=\overline{g_+^{-1}}\psi_-.
           \eeq
           If $k<0$, both sides of \eqref{5.5} must be equal to 0 because the left-hand side belongs to $\HH_s^+$ and vanishes at $i$,
           while the right-hand side is in  $\HH_s^-$  ($\frac{1}{s}=\frac{1}{p'}+\frac{1}{q}<1$).

           If $k\geq 0$, then both sides of \eqref{5.5} are equal to $\frac{P_k}{(\xi-i)^k}$, where $P_k$ denotes a polynomial of degree smaller or equal to $k$. Therefore we have
          $\psi_+=\frac{P_k}{(\xi+i)^k}\overline{g_-^{-1}}, \quad \psi_-=\frac{P_k}{(\xi-i)^k}\overline{g_+},$
           and  (ii) follows from here analogously as in the proof of Theorem \ref{thm5.1}.

\endpf

 Taking a simple example: let us denote $\left(\frac{\xi-i}{\xi+i}\right)^{1/2}_{\infty}=\sqrt{\frac{\xi-i}{\xi+i}}$. Then we have a canonical
$(q,2)$ - factorization, with $q>2$,
\beq \label{5.6}  \sqrt{\frac{\xi-i}{\xi+i}}=\sqrt{\xi-i}\,\, \frac{1}{\sqrt{\xi+i}}
\eeq
(with appropriate branches) where $g_-=\sqrt{\xi-i} \notin \HH_2^-$, $g_+^{-1}=\sqrt{\xi+i} \notin \HH_2^+,$  $g_-^{-1}=\frac{1}{\sqrt{\xi-i}}  \notin H_2^-,$ and
$g_+=\frac{1}{\sqrt{\xi+i}}  \notin H_2^+$. Therefore, by Theorems \ref{thm5.1} and \ref{thm5.2}, we have for the Toeplitz operator $T_{\sqrt{\frac{\xi-i}{\xi+i}}}$ in $H_2^+$
\beq \label{5.6A} \ker T_{\sqrt{\frac{\xi-i}{\xi+i}}}=\ker T^*_{\sqrt{\frac{\xi-i}{\xi+i}}}=\{0\}.
\eeq
Note that not only is \eqref{5.6} not a 2 - factorization, but moreover $\sqrt{\frac{\xi-i}{\xi+i}}$ does not admit any 2-factorization, whether or not of Wiener--Hopf type.
Note also that \eqref{5.6A} could also have been obtained from Theorem \ref{thm5.1} using the non-canonical $(q,2)$-factorization
$\sqrt{\frac{\xi-i}{\xi+i}}=\sqrt{\xi-i}\,\,\frac{\xi+i}{\xi-i}\,\,\frac{1}{\sqrt{\xi+i}}.$

Another example: consider the unbounded symbol $g(\xi)=\frac{\sqrt[6]{(\xi+i)^{14}}}{(\xi-i)^9}$. The function $g$ admits a (3,2)-factorization of the form
\[\frac{\sqrt[6]{(\xi+i)^{14}}}{(\xi-i)^9}=\sqrt{\xi-i}\left(\frac{\xi+i}{\xi-i}\right)^2\sqrt[3]{\xi+i}.\]
By Theorem 7.1, we have that $\dim \ker T_g=1$ in $H_2^+$ and $\ker T_g=\frac{1}{(\xi+i)^{4/3}}\spam \left\{\frac{1}{\xi+i}\right\}$.

% section8   {Products by non-integer powers of $r$ and piecewise continuous symbols
\section{Products by non-integer powers of $r$ and piecewise continuous symbols}
\label{sec:8}

We start by studying the relations between the kernels of two Toeplitz operators whose symbols  differ by a factor which is a non-integer power of $r$, motivated by the study of Toeplitz operators with piecewise continuous symbols.  For integer exponents,
these relations were studied in \cite{BCD, CMP}, in particular as regards their dimensions.

%Theor 8.1
\begin{thm} \label{thm8.1}
Let $f \in \sigma_p$ and let
\[g=fr_c^{\alpha} \quad \text{ with } c \in \RR, \, \, 0<\alpha <1\]
where $r_c^{\alpha}$ is defined in \eqref{r_form}. Then $D(T_g)=D(T_f)$ and
\[\ker T_g=\left(\frac{\xi+i}{\xi-c}\right)^{\alpha}\ker T_f \cap D(T_g).\]
\end{thm}
\beginpf
Let $\phi_+ \in H_p^+$. We have that
$\phi_+ \in \ker T_g \, \Leftrightarrow\, g\phi_+=\phi_- \in H_p^-$
and
\[g\phi_+=\phi_- \Leftrightarrow fr_c^{\alpha}\phi_+=\phi_- \Leftrightarrow f\underbrace{\left(\frac{\xi-c}{\xi+i}\right)^{\alpha}\phi_+}_{\in H_p^+}=\underbrace{\left(\frac{\xi-c}{\xi-i}\right)^{\alpha}\phi_-}_{\in H_p^-}\]
\[\Rightarrow \left(\frac{\xi-c}{\xi+i}\right)^{\alpha}\phi_+ \in \ker T_f \Rightarrow\phi_+ \in \left(\frac{\xi+i}{\xi-c}\right)^{\alpha}\ker T_f \cap D(T_g).\]
Conversely, let $\phi_+ \in \left(\frac{\xi+i}{\xi-c}\right)^{\alpha}\ker T_f \cap D(T_g)$. Then $\phi_+ =\left(\frac{\xi+i}{\xi-c}\right)^{\alpha}k_+$ where $k_+ \in \ker T_f$, with $fk_+=k_- \in H_p^-$, so
\[\underbrace{g\phi_+}_{\in L_p}=fr_c^{\alpha}\phi_+=f\left(\frac{\xi-i}{\xi-c}\right)^{\alpha}k_+=\left(\frac{\xi-i}{\xi-c}\right)^{\alpha}k_-=\underbrace{k_-/\left(\frac{\xi-c}{\xi-i}\right)^{\alpha}}_{\in \overline{\Nev_+}}\in L_p \cap\overline{\Nev_+}=H_p^-. \]
It follows that $\phi_+ \in \ker T_g$. \endpf

%Corol 8.2
\begin{cor} \label{cor8.2}
If $g=f\prod_{j=1}^m r_{c_j}^{\beta_j}$, where $c_1,c_2, \ldots, c_m$ are distinct points of $\RR$ and $0 <\beta_j <1 $ for all $j=1,2, \ldots, m$, then
\[\ker T_g = \prod_{j=1}^m\left(\frac{\xi+i}{\xi-c_j}\right)^{\beta_j}\ker T_f \cap D(T_g).\]
\end{cor}
Recall from Section 2 that piecewise continuous functions can be represented as products of a continuous function in $\RR_{\infty}$ by non-integer powers of $r$ and that, unlike continuous functions, a non-vanishing piecewise continuous function may not admit a Wiener--Hopf $p$-factorization, i.e., may be $p$-singular. In that case the range of the corresponding Toeplitz operator is not closed, so the operator is not Fredholm; we may ask however whether another type of representation of the symbol would allows us to describe the kernel of Toeplitz operator and its adjoint.

Indeed we show in the following theorem that for every symbol that can be represented in the form \eqref{g_decomp_b} where $h$ admits a bounded factorization, which includes in particular all non-vanishing piecewise H\"older continuous symbols, we can describe the corresponding Toeplitz kernels based on a
$(q,p')$-factorization (assuming the Toeplitz operator defined on $H_p^+$). Here we use the notation $\prod_{j=1}^0x_j=1$.

%Theor 8.3
\begin{thm} \label{thm8.3} Let $h=h_-r^kh_+$ with $h_{\pm} \in \mathcal{G}H^{\pm}_{\infty}, \, k \in \ZZ$, and
 \[g = h\, r_{\infty}^{\alpha_{\infty}}\left(\prod_{j=1}^{m}r^{\alpha_j}_{c_j}\right)\left(\prod_{j=1}^{s}r^{1/p}_{d_j}\right)\]
where $m,s$ are non-negative integer numbers, $c_j \, (j=1, \ldots, m)$ and $d_j \, (j=1, \ldots, s)$ are distinct real numbers,
\[-\frac{1}{p} < \alpha_{\infty} \leq \frac{1}{p'}, \quad -\frac{1}{p'} < \alpha_{j} < \frac{1}{p} \, \text{ for all } j=1, \ldots, m.\]
Then

\begin{description}
  \item[(i)] if $-\frac{1}{p} < \alpha_{\infty} < \frac{1}{p'}$, we have $\quad \dim \ker T_{g}= \max \left\{0, -k - s\right\}=\max \left\{0, \dim \ker T_h - s\right\}$ and, if $k < -s$,
      \[\ker T_g= h_+^{-1}\,\prod_{j=1}^s\left(\frac{\xi+i}{\xi-d_j}\right)^{1/p}\prod_{j=1}^m\left(\frac{\xi+i}{\xi-c_j}\right)^{\alpha_j}(\xi+i)^{\alpha_\infty}\,K_0\]
  with $K_0=\{\psi_+ \in K_{r^{|k|}}: \psi_+(d_j)=0, \, j=1,2, \ldots, s\}$.

  \item[(ii)] if $ \alpha_{\infty} = \frac{1}{p'}$, we have $\quad \dim \ker T_{g}= \max \left\{0, -k-s-1 \right\}=\max \left\{0, \dim \ker T_h-s-1 \right\}$ and, if $k<-s-1$,
 \[\ker T_g= h_+^{-1}\,\prod_{j=1}^s\left(\frac{\xi+i}{\xi-d_j}\right)^{1/p}\prod_{j=1}^m\left(\frac{\xi+i}{\xi-c_j}\right)^{\alpha_j}(\xi+i)^{1/p'}\,\widetilde{K_0}\]
  with $\widetilde{K_0}=\{\psi_+ \in K_{r^{|k|-1}}: \psi_+(d_j)=0, \, j=1,2, \ldots, s\}$.
\end{description}
\end{thm}
We prove this theorem in several steps, using the following Lemmas.

%Lemma 8.4
\begin{lem} \label{lem8.4} With the same assumptions as in Corollary \ref{cor8.2}, the function
\[f = h \,r_{\infty}^{\alpha_{\infty}}\left(\prod_{j=1}^m r_{c_j}^{\alpha_j}\right)\]
admits a $(q,p')$-factorization with $q \geq p$, of the form $f=f_-r^kf_+$ with
\[f_-=h_-(\xi-i)^{\alpha_{\infty}}\prod_{j=1}^m\left(\frac{\xi-i}{\xi-c_j}\right)^{\alpha_j}, \quad f_+=h_+\left(\frac{1}{\xi+i}\right)^{\alpha_{\infty}}\prod_{j=1}^m\left(\frac{\xi-c_j}{\xi+i}\right)^{\alpha_j}.\]
The function admits a $p$-factorization if and only if $\alpha_{\infty}\neq \frac{1}{p'}$.
\end{lem}
\beginpf We have
\[f=f_-r^kf_+\]
with
\[f_-=h_-(\xi-i)^{\alpha_{\infty}}\prod_{j=1}^{m}\left(\frac{\xi-i}{\xi-c_j}\right)^{\beta_j}\,\,,\,\,
f_+=h_+\frac{1}{(\xi+i)^{\alpha_{\infty}}}\prod_{j=1}^{m}\left(\frac{\xi-c_j}{\xi+i}\right)^{\beta_j}.\]
It is clear that $f_+ \in \HH_{p'}^+$, $f_-^{-1} \in \HH_{p'}^-$; on the other hand $f_- \in \HH_p^-$, $f_+^{-1} \in \HH_{p}^+$ if $-\frac{1}{p}< \alpha_{\infty} <\frac{1}{p'}$. If $\alpha_{\infty}=\frac{1}{p'}$ we have
\[\frac{f_-}{\xi-i}=\frac{h_-}{(\xi-i)^{1/p}}\prod_{j=1}^{m}\left(\frac{\xi-i}{\xi-c_j}\right)^{\beta_j}\]
with $-\frac{1}{p'}< \beta_j<\frac{1}{p}$ for all $j$.

Let $p<q<\min \{\frac{1}{\beta_j}, j=1, \cdots, m\}\cap \mathbb R^+$; then $\frac{f_-}{\xi-i}\in H_q^-$. Analogously, $f_+^{-1}\in \HH_q^+$.

It is left to prove that if $\alpha_{\infty}=\frac{1}{p'}$ then $f$ does not have a $p$-factorization (note that the non-existence of a WH $p$-factorization is a well known result). By Proposition \ref{prop3.6} it is enough to show that we cannot have
$f_-f_+=G_-rG_+$,
where the right-hand side is a $p$-factorization. Indeed we would have
\[f_-^{-1}=A\frac{G_-^{-1}}{\xi-i} \in H_{p'}^- \quad \text{with} \quad f_-^{-1}=h_-^{-1}\frac{1}{(\xi-i)^{1/p'}}\prod_{j=1}^{m}\left(\frac{\xi-c_j}{\xi-i}\right)^{\beta_j}\]
\[f_+=A\frac{G_+}{\xi+i} \in H_{p'}^+ \quad \text{with} \quad f_+=h_+\frac{1}{(\xi+i)^{1/p'}}\prod_{j=1}^{m}\left(\frac{\xi+i}{\xi-c_j}\right)^{\beta_j}\]
which is impossible because $h_+^{\pm 1} \in H_{\infty}^+, \, \, h_-^{\pm 1} \in H_{\infty}^-$. \endpf
Note that, in the previous result, $f_+^{-1}\left(\frac{\xi+i}{\xi-d_j}\right)^{1/p} \notin \HH_p^+$ for any $d_j \in \RR$.

%Lemma 8.5
\begin{lem} \label{lem8.5} Let $d_j$, $j=1,2, \ldots, s$, be distinct points in $\RR$ and let
$g= f \left(\prod_{j=1}^s r^{1/p}_{d_j}\right)$,
where we assume that $f$ admits a $(q, p')$-factorization with $q \geq p$, of the form
\beq \label{eq1_lem8.5}
f=f_-r^kf_+
\eeq
such that, for all $j=1, 2, \ldots, s$, we have
%\beq \label{eq2_lem8.5}
$f_+^{-1}\left(\frac{\xi+i}{\xi-d_j}\right)^{1/p} \notin \HH_p^+$.
%\eeq
Then
\begin{description}
  \item[(i)] $ \ker T_g=0$ if $s \geq \dim \ker T_f$,
  \item[(ii)] $ \ker T_g= \prod_{j=1}^s\left(\frac{\xi+i}{\xi-d_j}\right)^{1/p}f_+^{-1}K_0$ with $K_0=\{\psi_+ \in K_{r^{|k|}}: \psi_+(d_j)=0, \, j=1,2, \ldots, s\}$ if $s < \dim \ker T_f$ and \eqref{eq1_lem8.5} is a $p$-factorization,
  \item[(iii)] $ \ker T_g= \prod_{j=1}^s\left(\frac{\xi+i}{\xi-d_j}\right)^{1/p}\frac{f_+^{-1}}{\xi+i}\widetilde{K_0}$ with $\widetilde{K_0}=\{\psi_+ \in K_{r^{|k|-1}}: \psi_+(d_j)=0, \, j=1,2, \ldots, s\}$ if $s < \dim \ker T_f$ and either $f_- \notin \HH_p^-$ or $f_+^{-1} \notin \HH_p^+$.
      \item[(iv)] $\dim \ker T_g = \max \{0, \dim \ker T_f -s\}$.
\end{description}
\end{lem}
\beginpf From Corollary \ref{cor8.2},
\[\ker T_g =\prod_{j=1}^s\left(\frac{\xi+i}{\xi-d_j}\right)^{1/p}\ker T_f \cap H_p^+.\]
If \eqref{eq1_lem8.5} is a $p$-factorization then, from Theorem \ref{thm5.1} and Corollary \ref{cor8.2}
\[\ker T_g=\{0\}, \, \text{ if } k \geq 0,\]
\[\ker T_g= \prod_{j=1}^s \left( \frac{\xi+i}{\xi-d_j} \right)^{1/p} f_+^{-1} K_{r^{|k|}}\cap H_p^+, \, \text{ if } k < 0.\]
So, for $k<0$, since $f_+^{-1} \in \HH_p^+$ and $\left(\frac{\xi+i}{\xi-d_j}\right)^{1/p}f_+^{-1} \notin \HH_p^+$, we see that
\[\ker T_g= \prod_{j=1}^s\left(\frac{\xi+i}{\xi-d_j}\right)^{1/p}f_+^{-1}K_0 \]
where $K_0=\{\psi_+ \in K_{r^{|k|}}: \psi_+(d_j)=0, \, j=1,2, \ldots, s\}$ is equal to $\{0\}$ if $s \geq |k|=-k$ and $\dim K_0=|k|-s$ if $s < |k|$. Thus we have $\ker T_g=\{0\}$ if $s \geq -k$, $\dim \ker T_g=|k|-s$ if $s<-k$, where $-k=|k|=\dim \ker T_f$ by Theorem \ref{thm5.1} if $k<0$.

If \eqref{eq1_lem8.5} is not a $p$-factorization, i.e., $f_- \notin \HH_p^-$ or $f_+^{-1} \notin \HH_p^+$, then we conclude analogously from Theorem \ref{thm5.1} and Corollary \ref{cor8.2} that
\[\ker T_g=\{0\}, \, \text{ if } k \geq -1,\]
\[\ker T_g= \prod_{j=1}^s\left(\frac{\xi+i}{\xi-d_j}\right)^{1/p}\frac{f_+^{-1}}{\xi+i}K_{r^{|k|-1}}\cap H_p^+=\prod_{j=1}^s\left(\frac{\xi+i}{\xi-d_j}\right)^{1/p}\frac{f_+^{-1}}{\xi+i}\widetilde{K_0}\]
where $\widetilde{K_0}=\{\psi_+ \in K_{r^{|k|-1}}: \psi_+(d_j)=0, \, j=1,2, \ldots, s\}$, if $k<-1$.

Thus $\ker T_g=\{0\}$ if $s \geq -k-1$ and $\dim \ker T_g=|k|-1-s$ if $s<-k-1$, where $-k-1=|k|-1=\dim \ker T_f$ by Theorem \ref{thm5.1} if $k<-1$.
\endpf

\beginpf \emph{(of Theorem \ref{thm8.3})} Let $f=h\, r_{\infty}^{\alpha_{\infty}}\left(\prod_{j=1}^m r_{c_j}^{\alpha_j}\right)$. By Lemma \ref{lem8.4}, $f$ has a $(q,p')$-factorization with $q \geq p,$ which is a $p$-factorization if and only if $\alpha_{\infty}\neq \frac{1}{p'}$. Thus, by Lemma \ref{lem8.5}, % since in this case \eqref{eq2_lem8.5} holds,
(i) and (ii) hold. \endpf

The result of Theorem \ref{thm8.3} has a simple geometric interpretation in the case of piecewise continuous symbols. Given any piecewise continuous function in $\RR_{\infty}$ of the form \eqref{g_decomp}, with discontinuity at the points $c_j$, $j=1,2, \ldots, n$, and (possibly) at $\infty$, the image of the continuous factor $h$ in the complex plane is a closed contour which does not pass by the origin, with winding number $k$, while the image of the function $g_p$ associated to $g$, given by \eqref{g_p}, is a closed curve that includes the image of $g$ as well as arcs connecting $g(\xi^-)$ and $g(\xi^+)$ whenever these are different. If a point of discontinuity $\xi=c_j \in \RR$ is such that the corresponding exponent $\alpha_j$ in \eqref{g_decomp} is $1/p$ or, for the point $\xi=\infty$, if the corresponding exponent $\alpha_{\infty}$ is $1/p'$, then the curve connecting $g(\xi^-)$ and $g(\xi^+)$ passes by the origin.

Thus we can interpret the result of Theorem \ref{thm8.3} as saying that $\dim \ker T_g$ is obtained from the dimension of $\ker T_h$ (which is zero if $k\geq 0$, and equal to $|k|$ if $k<0$) in the following way: if $k \geq 0,$ $\dim \ker T_g=0$; if $k<0$, then $\dim \ker T_g$ is obtained by subtracting from $|k|$ the number of arcs passing by the origin in the image of $g_p$, if the number of these arcs is smaller than $|k|$; otherwise $\dim \ker T_g=0.$

We present below some examples illustrating this geometric interpretation.

\begin{enumerate}
  \item Consider the function $g$, given by
  \[
  g(\xi)=\frac{(\xi+i)^3}{(\xi-2i)(\xi-3i)(\xi-4i)}\sqrt{\frac{\xi-i}{\xi+i}},
  \]
  which is a piecewise continuous function with a discontinuity at $c=\infty$.  We plot the image of $g_2$ (see \eqref{g_2})  in the complex plane, which includes   $\{g(\xi): \xi \in \RR\}$ as well as the segment connecting $g(\infty^+)$ and $g(\infty^-)$, which  passes by the origin.  By Theorem \ref{thm8.3},
$\dim \ker T_g=2$ in $H_2^+$, since $\dim \ker T_h=3$ where
$
h(\xi)=\frac{(\xi+i)^3}{(\xi-2i)(\xi-3i)(\xi-4i)},
$
 and we have
 %\[\ker T_g=\left(\frac{\xi}{\xi+i}\right)^{1/3}\widetilde{K_0},\]
$\ker T_g=(\xi+i)^{1/2}K_{r^{2}}.$
%with $\widetilde{K_0}=\left\{\psi_+ \in K_{r^{2}}: \psi_+(0)=0\right\}$.

\begin{figure}[H]
\centering
\includegraphics[width=.8\linewidth]{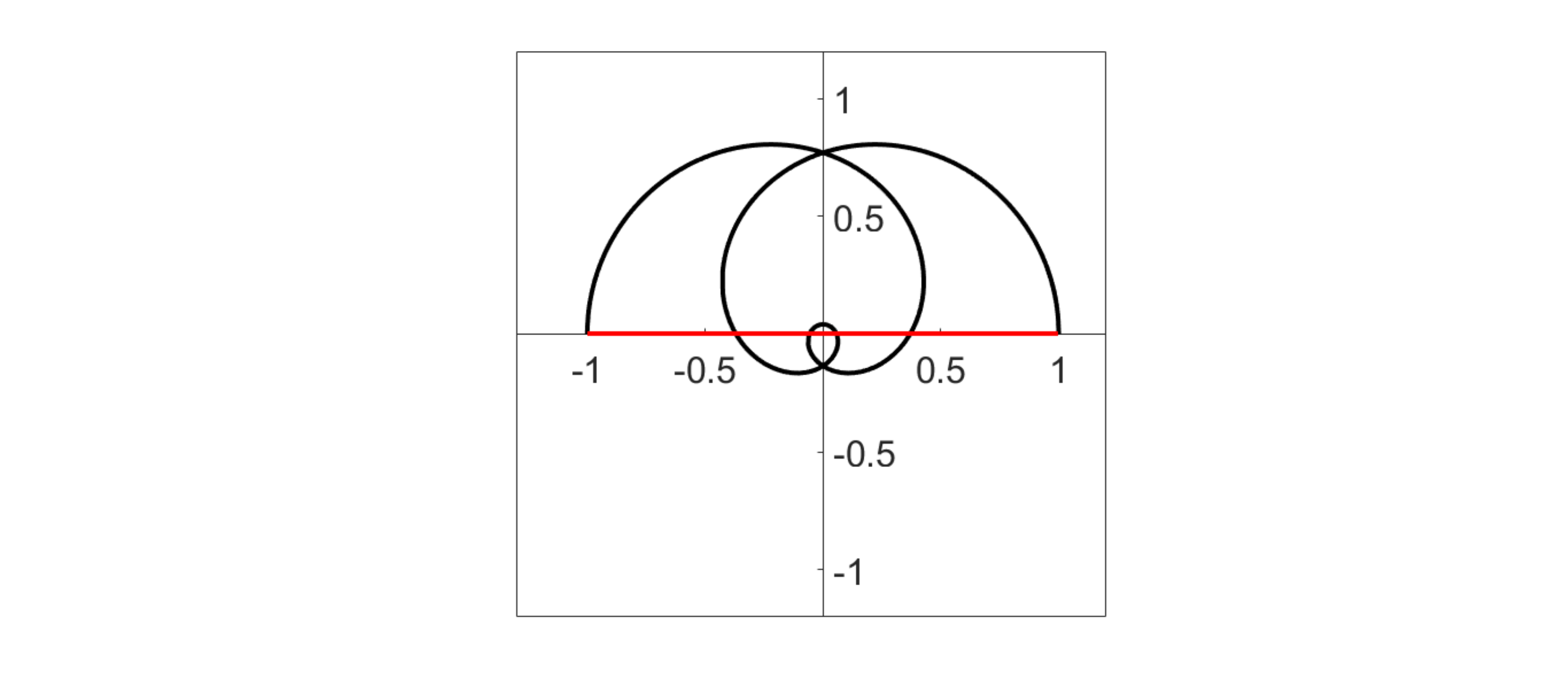}
\caption{Plot of $g_2$.}
\label{figures}
\end{figure}

\item Let now $G$ be given by
\[
G(\xi)=\sqrt{\frac{\xi-i}{\xi+i}}\left(\frac{(\xi+i)^2}{(\xi-2i)(\xi-3i)}\right)\left(\frac{\xi-i}{\xi+i}\right)^{1/2}_0. \] $G$ is a piecewise continuous function with a discontinuity at $c=0$ and $c=\infty$; the image of  $G_2$, obtained according to \eqref{g_2}, is shown in the figure below.
%By Theorem %\ref{thm8.3}
We have that $\dim \ker T_G=0$ in $H_2^+$ since $\dim \ker T_h=2$, where $h=\frac{(\xi+i)^2}{(\xi-2i)(\xi-3i)}$.

\begin{figure}[H]
\centering
\includegraphics[width=.8\linewidth]{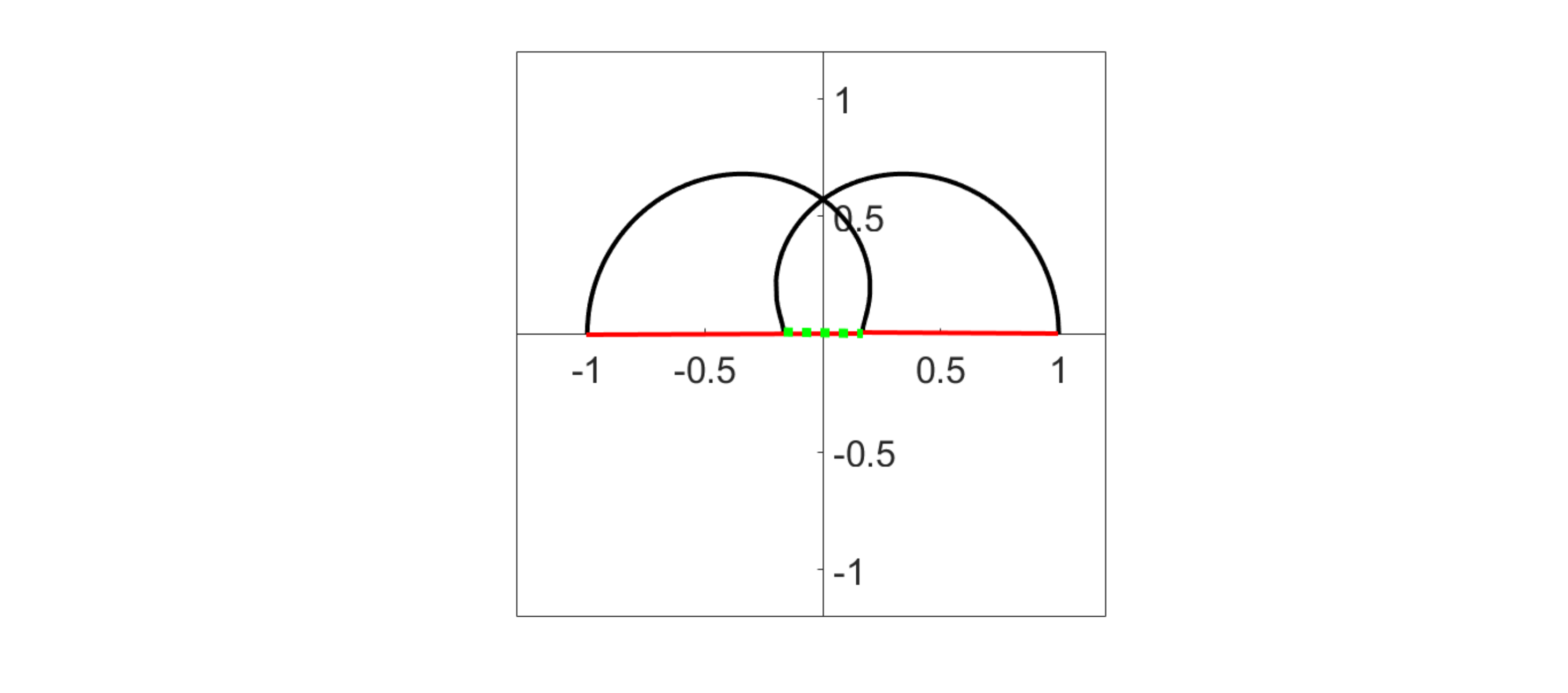}
\caption{Plot  of $G_2$.}
\label{figures}
\end{figure}

\end{enumerate}

\subsection*{Acknowledgements}
This work
was partially supported by FCT/Portugal through UID/MAT/04459/2020.
 The research of  M. T. Malheiro  was partially supported by Portuguese
Funds through FCT/Portugal within
the  Projects UIDB/00013/2020
and UIDP/00013/2020.

\end{document}